\documentclass[preprint,10pt]{elsarticle}

\usepackage{amsmath,amssymb,amsfonts,amsthm}
\usepackage{graphics,graphicx,epsfig,float}
\usepackage{psfrag}
\usepackage{cases}
\usepackage[english]{babel}
\usepackage[T1]{fontenc}
\usepackage{textcomp}
\usepackage{ae}
\usepackage{aecompl}
\usepackage[cm]{aeguill}
\usepackage{color}
\usepackage{array}

\newcommand{\vw}{{\mathbf{w}}}

\newcommand{\vf}{{\mathbf{f}}}

\newcommand{\vS}{{\mathbf{S}}}

\newcommand{\cT}{{\mathcal T}}
\newcommand{\cQ}{{\mathcal Q}}
\newcommand{\cR}{{\mathcal R}}

\newcommand{\R}{\ensuremath{\mathbb{R}}}

\newcommand{\dt}{\partial_t}
\newcommand{\dx}{\partial_x}
\newcommand{\eps}{\varepsilon}
\newcommand{\dsp}{\displaystyle}

\newcommand{\beq}{\begin{equation}}
\newcommand{\eeq}{\end{equation}}
\newcommand{\be}{\begin{equation}}
\newcommand{\ee}{\end{equation}}
\newcommand{\bes}{\begin{equation*}}
\newcommand{\ees}{\end{equation*}}
\newcommand{\bea}{\begin{eqnarray}}
\newcommand{\eea}{\end{eqnarray}}
\newcommand{\beas}{\begin{eqnarray*}}
\newcommand{\eeas}{\end{eqnarray*}}
\newcommand{\beqn}{\begin{eqnarray}}
\newcommand{\eeqn}{\end{eqnarray}}
\newcommand{\beqns}{\begin{eqnarray*}}
\newcommand{\eeqns}{\end{eqnarray*}}

\def\fr#1#2{\frac{#1}{#2}}

\newtheorem{thrm}{Theorem}[section]

\newtheorem{rmrk}[thrm]{Remark}

\newtheorem{nota}[thrm]{Notation}
\newtheorem{proposition}{Proposition}

\begin{document}

\begin{frontmatter}

\title{A splitting approach for the fully nonlinear and weakly dispersive Green-Naghdi model}

\author[label1]{P. Bonneton}
\address[label1]{Universit\'e Bordeaux 1, CNRS, UMR 5805-EPOC, F-33405 Talence, France}
\ead{p.bonneton@epoc.u-bordeaux1.fr}
\author[label2]{F. Chazel}
\address[label2]{Universit\'e de Toulouse, UPS/INSA, IMT, CNRS UMR 5219, F-31077 Toulouse, France}
\ead{florent.chazel@math.univ-toulouse.fr}
\author[label3]{D. Lannes}
\address[label3]{DMA, Ecole Normale Sup\'erieure et CNRS UMR 8553, 45 rue d'Ulm, F-75005 Paris, France}
\ead{david.lannes@ens.fr}
\author[label4]{F. Marche}
\address[label4]{I3M, Universit\'e Montpellier 2, CC 051, F-34000 Montpellier, France}
\ead{fabien.marche@math.univ-montp2.fr}
\author[label1]{M. Tissier}
\ead{m.tissier@epoc.u-bordeaux1.fr}

\journal{J. Comp. Phys.}

\begin{abstract}
The fully nonlinear and weakly dispersive Green-Naghdi model for shallow water waves of large amplitude is studied. The original model is first recast under a new formulation more suitable for numerical resolution. An hybrid finite volume and finite difference splitting approach is then proposed. The hyperbolic part of the equations is handled with a high-order finite volume scheme allowing for breaking waves and dry areas. The dispersive part is treated with a classical finite difference approach. Extensive numerical validations are then performed in one horizontal dimension, relying both on analytical solutions and experimental data. The results show that our approach gives a good account of all the processes of wave transformation in coastal areas: shoaling, wave breaking and run-up.
\end{abstract}


\end{frontmatter}

\section{Introduction}\label{intro}

In an incompressible, homogeneous, inviscid fluid, the propagation of surface waves is governed by the 
Euler equations with nonlinear boundary conditions at the surface and at the bottom. In its full generality, 
this problem is very complicated to solve, both mathematically and numerically. This is the reason why more simple 
models have been derived to describe the behavior of the solution in some physical specific regimes. A recent review of 
the different models that can be derived can be found in \cite{Lannes-Bonneton}.\\
Of particular interest in coastal oceanography is the \emph{shallow-water} regime, which corresponds to the configuration
where the wave length $\lambda$ of the flow is large compared to the typical depth $h_0$: 
$$
\mbox{(Shallow water regime) }\quad \mu:=\frac{h_0^2}{\lambda^2}\ll 1.
$$ 
When the typical amplitude $a$ of the wave is small, in the sense that 
$$
\mbox{(Small amplitude regime) }\quad \eps:=\frac{a}{h_0}=O(\mu),
$$
 it is known that an approximation of order $O(\mu^2)$ of the free surface Euler equations is furnished by the Boussinesq systems, such as the one derived
by Peregrine in \cite{per1967} for uneven bottoms. This model couples the surface elevation $\zeta$ to the vertically averaged horizontal component of the velocity $V$, and can be written in non-dimensionalized form as
\begin{equation}\label{Bo1}
\left\lbrace
\begin{array}{l}
\partial_t\zeta+\nabla\cdot (h V)=0,\\
\partial_t V+\eps (V\cdot \nabla)V+\nabla\zeta=\mu {\mathcal D}+O(\mu^2),
\end{array}\right.
\end{equation}
where $h$ is the water depth and ${\mathcal D}$ accounts for the nonhydrostatic and dispersive effects, and is a function of $\zeta$, $V$ and their
derivatives. For instance, in the Boussinesq model derived in \cite{per1967}, one has
\begin{equation}\label{D}
{\mathcal D}=\frac{h}{2}\nabla [\nabla\cdot(h\dt V)]-\frac{h^2}{6}\nabla^2\dt V.
\end{equation}
Unfortunately, the small amplitude assumption $\eps=O(\mu)$ is too restrictive for many applications in coastal oceanography,
where \emph{large amplitude} waves have to be considered,
$$
\mbox{(Large amplitude regime) }\quad \eps:=\frac{a}{h_0}=O(1).
$$
If one wants to keep the same $O(\mu^2)$ precision of (\ref{Bo1}) in a large amplitude regime, then the expression for
${\mathcal D}$ is much more complicated than in (\ref{D}). For instance, in $1D$ and for flat bottoms, one has
$$
{\mathcal D}=\frac{\eps}{3h}\dx[h^3(V_{xt}+VV_{xx}-(V_x)^2)].
$$
In this regime, the corresponding equations (\ref{Bo1}) have been derived first by Serre and then Su and Gardner \cite{sug1969}, Seabra-Santos \emph{et al.} \cite{sea1987} and Green and Naghdi \cite{gre1976} (other relevant references are \cite{cie2006,wei1995,MS}); consequently, these equations carry several names: Serre, Green-Naghdi, or fully nonlinear Boussinesq equations.  We will call them \emph{Green-Naghdi equations} throughout this paper. Here again, we refer to \cite{Lannes-Bonneton} for more details; note also that a rigorous mathematical justification of these models has been given in \cite{AL}.

\medbreak

The Green-Naghdi equations (\ref{Bo1}) provide a correct description of the waves up to the breaking point; from this point 
however, they become useless (at least without consequent modifications). A first approach to model wave breaking 
is to add an ad hoc viscous term to the momentum equation, whose role is to 
account for the energy dissipation that occurs during wave breaking. This ap- 
proach has been used for instance by Zelt \cite{Zelt} or Kennedy \cite{Kennedy} and Chen \cite{Chen}. Recently, 
Cienfuegos et al. \cite{cie2009} proposed a new 1D wave-breaking parametrization 
including viscous-like effects on both the mass and the momentum equations.
This approach is able to reproduce wave height decay and intraphase nonlinear 
properties within the entire surf zone. However, the extension of this ad hoc 
parametrization to 2D wave cases remains a very difficult task.
Another approach to handle wave breaking is to use the classical nonlinear shallow water equations, defined with ${\mathcal D}=0$ in (\ref{Bo1}) and denoted by NSWE in the following. 
These equations being hyperbolic, they develop shocks; after the breaking point, the waves are then described 
by the \emph{weak solutions} of this hyperbolic system. This approach, used in \cite{koba1989}  and \cite{bonn2007} is satisfactory in the sense that it gives a natural and correct description of the dissipation of energy during wave 
breaking. Its drawback, however, is that it is inappropriate
 in the shoaling zone since this models neglects the nonhydrostatic and dispersive effects. The motivation of 
this paper is to develop a model and a numerical scheme that describes correctly both phenomena. More precisely, we want to
\begin{enumerate}
\item[$\sharp$ 1] Provide a good description of the dispersive effects (in the shoaling zone in particular);
\item[$\sharp$ 2] Take into account wave breaking in a simple way.
\end{enumerate}

\medbreak

Another theoretical and numerical difficulty in coastal oceanography is the description of the shoreline, i.e. the zone where the water depth vanishes, as the size of the computational domain becomes part of the solution. Taking into account the possibility of a vanishing depth while keeping the dispersive effects is more difficult, see \cite{cie2007,funwave} for instance.
As for the breaking of waves, neglecting the nonhydrostatic and dispersive
effects makes the things simpler. Indeed, various efficient schemes have been developed to handle the possibility of vanishing depth for the NSWE with source terms, relying for instance on coordinates transformations \cite{brocchini}, artificial porosity \cite{vanhof} or even variables extrapolations \cite{lynett}. In a simpler way, it is shown in  \cite{marche_bonneton} that the occurrence of dry areas can be naturally handled with a water height positivity preserving finite volume scheme, without introducing any numerical trick. Of course, the price to pay is the same as above: the dispersive effects are lost. Hence the third motivation for this paper:
\begin{enumerate}
\item[$\sharp$ 3] Propose a simple numerical method that allows at the same time the possibility of vanishing depth and dispersive effects.
\end{enumerate}

\medbreak

The strategy adopted here to handle correctly the three difficulties $\sharp$ 1-3 
identified above starts from the Green-Naghdi equations. 
As already said, they are very well adapted to $\sharp$ 1. With a careful choice of the numerical methods, they also
allow for the possibility of vanishing depth, and thus answer to $\sharp$ 3. The main difficulty is thus to handle $\sharp$ 2 (i.e. wave breaking) with a code based on the Green-Naghdi equation. In order to do so, we use a numerical scheme that decomposes the hyperbolic and dispersive parts of the equations \cite{erduran}. We also refer to \cite{LGH} for a recent  numerical analysis
of the Green-Naghdi equations based on a Godunov type scheme and that provide good results
for the dam break problem. We use here a second order splitting scheme, we compute the approximation $U^{n+1}=(\zeta^{n+1},V^{n+1})$ at time $(n+1)\delta_t$ in
terms of the approximation $U^n$ at time $n\delta_t$ by solving
$$
U^{n+1}=S_1(\delta_t/2)S_2(\delta_t)S_1(\delta_t/2)U^n,
$$
where $S_1(\cdot)$ is the solution operator associated to the NSWE (${\mathcal D}=0$ in (\ref{Bo1})) and
$S_2(\cdot)$ the solution operator associated to the dispersive part of the equations (keeping only the time derivatives and the r.h.s. of (\ref{Bo1})). For the numerical computation of $S_1(\cdot)$, we use a high order, robust and well-balanced finite volume method, based on a relaxation approach \cite{berthon_marche}. This method is known to be computationally cheap and very efficient to handle wave breaking and presents another interesting feature for our purposes: it allows the localisation of the shocks. In the vicinity of these shocks (or bores to use the physical term) the derivation of the dispersive components of the Green-Naghdi equation is meaningless and these terms, which contain third order derivatives, become very singular; moreover, it is known \cite{bonn2007,brocchini2} that the NSWE correctly describe the dynamics of the waves near the breaking point. We therefore ``skip'' the computation of $S_2(\cdot)$ near the shocks detected during the computation of $S_1(\delta_t/2)$. Elsewhere, $S_2(\cdot)$ is computed using a finite difference scheme (note that a careful mathematical analysis of $S_2(\cdot)$ allows considerable simplifications and numerical improvements).

\medbreak

In Section \ref{model}, we present the physical model studied here, namely, the Green-Naghdi equations. After giving
the formulation of the equations in non-dimen\-sionalized form in \S \ref{NDGN}, we show in \S \ref{reform} that it is possible
to rewrite them in a convenient way that does not require the computation of any third order derivative of the unknowns $\zeta$ and $V$, and exploits the regularizing
properties of the equations. With classical methods, we then turn to derive a family of Green-Naghdi equations with improved frequency dispersion in \S \ref{sectimp}, depending
on a parameter $\alpha$ to be chosen. Another, still non-dimensionlized, reformulation of the equations (in terms of $(h,hV)$ rather than $(\zeta,V)$) is then
given in \S \ref{secthhV}; finally, we give a version with dimensions of these equations in \S \ref{sectdim}.

Section \ref{sectNS} is then devoted to the presentation of the numerical scheme. The hyperbolic/dispersive splitting
is introduced in \S \ref{sectsplit}. We then turn to describe the spacial discretization of the hyperbolic and dispersive parts in \S \ref{secthyp} and \S \ref{sectdisp} respectively. The time discretization is described in \S \ref{secttime}, where the
consequences of our approach for the dispersive properties of the model are also studied carefully. We show in particular that it is possible to derive an exact formula for the
semi-discrete dispersion relation that approaches the exact one at order two. We then
use this formula to choose the best parameter $\alpha$ for the
frequency-improved Green-Naghdi equations; this choice
differs from the classical one based on the exact dispersion relation.

Finally, we present in Section \ref{numval} several numerical validations of our model. We first consider the case of solitary waves in \S \ref{solwav} and use it as a validation tool for our numerical scheme. We then evaluate the dispersive properties of the model by considering the propagation of a periodic and regular wave over a flat bottom in \S \ref{fenton};
this test illustrates the interest of choosing the frequency parameter $\alpha$
in terms of the semi-discrete dispersion relation. In \S \ref{wallwb}, we focus on the reflection of a solitary wave at a wall, while the ability of the model to simulate the nonlinear shoaling of solitary waves over regular sloping beaches is investigated in \S \ref{LEGI}. At last, the run-up and run-down of a breaking solitary wave is studied in \S \ref{syno}.

\section{The physical model}\label{model}

\begin{figure}[t]
\psfrag{z}{{\footnotesize $z$}}
\psfrag{0}{{\footnotesize $0$}}
\psfrag{h0}[c][c]{\textcolor{white}{7\footnotesize{2}}{\footnotesize $-h_0$}}
\psfrag{N}{{\footnotesize $z=\zeta(t,X)$}}
\psfrag{F}{{\footnotesize $z=-h_0+b(X)$}}
\psfrag{X}{{\footnotesize $X=(x,y)$}}
\hspace{0.45em}\includegraphics[width=0.85\textwidth]{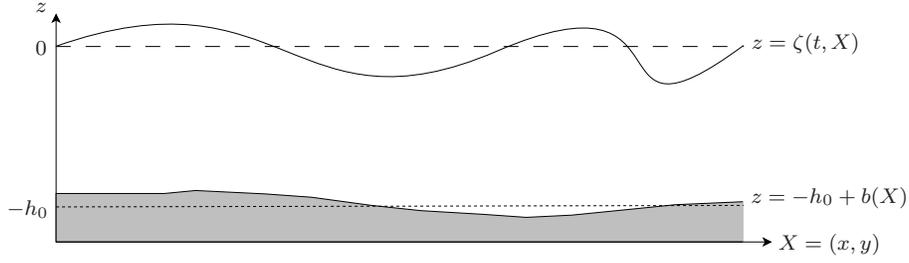}
\caption{Sketch of the domain}
\label{domain}
\end{figure}

Throughout this paper, we denote by $\zeta(t,X)$ the elevation of the surface with respect to its rest state, and by $-h_0+b(X)$ a parametrization of the bottom, where $h_0$ is a reference depth (see Figure \ref{domain}). Here $X$ stands for the horizontal variables ($X=(x,y)$ for $2D$ surface waves, and $X=x$ for $1D$ surface waves), and $t$ is the time variable; we also denote by $z$ the
vertical variable.

\noindent If $U_{hor}$ denotes the horizontal component of the velocity field in the fluid domain, we then define $V$ as
$$
V(t,X)=\frac{1}{h}\int_{-h_0+b}^\zeta U_{hor}(t,X,z)dz,
$$
where $h:=h_0+\zeta-b$ is the water depth. We thus have $V=(u,v)\in \R^2$ for $2D$ surface waves, and $V=u\in\R$ for $1D$ surface waves.

\smallbreak

\noindent Denoting by $a$ the typical amplitude of the waves, by $a_{bott}$ the typical amplitude of the bottom variations,
and by $\lambda$ the order of the wavelength of the wave, it is possible to define
dimensionless variables and unknowns as
$$
\widetilde X=\frac{X}{\lambda},\qquad \widetilde t= \frac{\sqrt{gh_0}}{\lambda} t
$$
and
$$
\widetilde \zeta=\frac{\zeta}{a},\qquad \widetilde b=\frac{b}{a_{bott}},\qquad \widetilde V=\frac{V}{\sqrt{gh_0}}.
$$
We also define three dimensionless parameters as
$$
\eps=\frac{a}{h_0},\qquad \mu=\frac{h_0^2}{\lambda^2},\qquad \beta=\frac{a_{bott}}{h_0};
$$
here $\eps$ denotes the nonlinearity parameter, $\mu$ is the shallowness parameter while $\beta$ accounts for the topography variations.

\subsection{The non-dimensionalized Green-Naghdi equations}\label{NDGN}

According to \cite{AL,Lannes-Bonneton}, the Green-Naghdi equations can be written under the following non-dimensionalized form 
(we omit the tildes for dimensionless quantities for the sake of clarity):
\be\label{eq1}
	\left\lbrace
	\begin{array}{l}
	\dsp \dt \zeta +\nabla\cdot (h V)=0,\\
	\dsp (I+\mu\cT[h,b])\dt V+\nabla\zeta+\eps(V\cdot \nabla)V
	+\eps\mu \cQ[h,b](V)=0,
	\end{array}\right.
\ee
where we still denote by $h$ the non-dimensionalized water depth,
$$
	h=1+\eps\zeta-\beta b,
$$
and the linear operator
$\cT[h,b]\cdot$ and the quadratic form $\cQ[h,b](\cdot)$ 
are defined for all smooth enough $\R^d$-valued function $W$ ($d=1,2$ is the surface dimension) by
\bea
	\label{eq2}
		\!\!\!\!\!\!\cT[h,b]W\!\!\!\!&=&\!\!\!\!\!
	\cR_1[h,b](\nabla\cdot W)
	+\beta\cR_2[h,b](\nabla b\cdot W)\\
	\label{eq3}
	\!\!\!\!\!\!\cQ[h,b](W)\!\!\!\!&=&\!\!\!\!\!
	\cR_1[h,b](\nabla\cdot(W\nabla\cdot W)-2(\nabla\cdot W)^2)
	+\beta\cR_2[h,b]((W\cdot\nabla)^2 b),
\eea
with, for all smooth enough scalar-valued function $w$,
\bea
	\label{eq4}
	\cR_1[h,b]w&=&-\frac{1}{3h}\nabla(h^3 w)
	-\beta\frac{h}{2}w\nabla b,\\
	\label{eq5}
	\cR_2[h,b]w&=&\frac{1}{2h}\nabla(h^2 w)
	+\beta w\nabla b.
\eea

\begin{nota}
	For the sake of clarity, and one no confusion is possible,
	we often write $\cT$, 
	$\cQ$, $\cR_1$ and $\cR_2$
	instead of $\cT[h,b]$, 
	$\cQ[h,b]$, etc.
\end{nota}

\begin{rmrk}\label{friction}
For practical applications (see for instance in \S \ref{syno}), a classical quadratic friction term can be added to the right-hand side of the momentum equation. It has the following expression: $-f \eps \mu^{-1/2} \frac{1}{h} \|V\| V$, where $f$ is a non-dimensional friction coefficient. 
\end{rmrk}
\subsection{First reformulation of the equations}\label{reform}

Let us now remark that 
$$
	\cQ(V)=\cT((V\cdot\nabla)V)+\cQ_1(V),
$$
where $\cQ_1(V)$ only involves second order derivatives of $V$ (while third order
derivatives appear in $\cQ(V)$); indeed,  a close look at
(\ref{eq2}) and (\ref{eq3}) shows that
\begin{eqnarray*}
	\cQ_1(V)&=&
	\cR_1\big(\nabla\cdot(V\nabla\cdot V-(V\cdot\nabla)V)-2(\nabla\cdot V)^2\big)\\
	& &+\beta\cR_2(V\cdot(V\cdot\nabla)\nabla b)\\
	&=&-2\cR_1(\partial_1V\cdot\partial_2V^\perp+(\nabla\cdot V)^2)
	+\beta\cR_2(V\cdot(V\cdot\nabla)\nabla b),
\end{eqnarray*}
with $V^\perp=(-V_2,V_1)^T$. The fact that this expression does not involve third order derivatives
is of great interest for the numerical applications.\\
We have thus obtained the following equivalent formulation of the Green-Naghdi equations (\ref{eq1}):
\be\label{eq6}
	\left\lbrace
	\begin{array}{l}
	\dsp \dt \zeta +\nabla\cdot (h V)=0,\\
        \dsp (I+\mu\cT)\dt V+\eps 	(I+\mu\cT)(V\cdot \nabla)V+
	\nabla\zeta 
	+\eps\mu\cQ_1(V)=0,
	\end{array}\right.
\ee
with 
$
	h=1+\eps\zeta-\beta b, 
$
and where
 the quadratic form $\cQ_1$ is given by
\be\label{eq11}
	\cQ_1[h,b](V)=-2\cR_1(\partial_1V\cdot\partial_2V^\perp+(\nabla\cdot V)^2)
	+\beta\cR_2(V\cdot(V\cdot\nabla)\nabla b)
\ee
(the linear operators $\cT$, $\cR_1$ and $\cR_2$
being defined in (\ref{eq2}), (\ref{eq4}) and (\ref{eq5}).

\subsection{Green-Naghdi equations with improved dispersive properties}\label{sectimp}

It is classical \cite{wit1984,mad1991}  or \cite{cie2007} that the frequency dispersion of (\ref{eq6}) can be improved
by adding some terms of order $O(\mu^2)$ to the momentum equation. Since this equation is
already precise up to terms of order $O(\mu^2)$, this manipulation does not affect the 
precision of the model. Such a manipulation is also performed in \cite{LGH} but with
the goal of working with potential variables rather than improving the frequency dispersion.\\
The first step consists in noticing that, from the second equation in (\ref{eq6}), one has
$$
\partial_t V = -\nabla\zeta - \eps (V\cdot \nabla)V+O(\mu),
$$
and therefore, for any parameter $\alpha\in \R$,
$$
\partial_t V=\alpha\partial_t V-(1-\alpha)\big(\nabla\zeta+\eps(V\cdot\nabla)V\big) +O(\mu).
$$
Replacing $\partial_t V$ by this expression in (\ref{eq6}) and dropping the $O(\mu^2)$
terms yields the following Green-Naghdi equations with improved frequency dispersion,
\be\label{eq6imp}
	\left\lbrace
	\begin{array}{l}
	\dsp \dt \zeta +\nabla\cdot (h V)=0,\\
        \dsp (I+\mu\alpha\cT)\dt V+\eps (I+\mu\alpha\cT)(V\cdot \nabla)V\\
        \indent\indent\dsp +
	(I-\mu(1-\alpha)\cT)\nabla\zeta 
	+\eps\mu\cQ_1(V)=0.
	\end{array}\right.
\ee
Of course, (\ref{eq6}) corresponds to a particular case of (\ref{eq6imp}) with $\alpha=1$.
The interest of working with (\ref{eq6imp}) is that it allows to improve the dispersive properties of the model by minimizing - thanks to the parameter $\alpha$ - the phase velocity error (see \ref{sectdim}). In \cite{proc}, a three-parameter family of formally equivalent Green-Naghdi equations is derived yielding further improvements of the dispersive properties. For the sake of simplicity, we stick here to the one-parameter family of Green-Naghdi systems (\ref{eq6imp}).

\subsection{Reformulation in terms of the $(h,hV)$ variables}\label{secthhV}

The Green-Naghdi equations with improved dispersion (\ref{eq6imp}) are stated as two
evolution equations for $\zeta$ and $V$. It is possible to give an 
equivalent formulation as a system of two evolution equations on $h$ and $hV$, as shown in 
this section.\\
For the first equation, one just has to remark that $\eps \dt\zeta=\dt h$, so that
$$
\dt h+\eps \nabla\cdot (hV)=0.
$$
 For the second equation, we first use this identity  to remark that 
$h\dt V=\dt (hV)+\eps\nabla\cdot (hV)V$. Multiplying the 
second equation of (\ref{eq6imp}) by $h$, and using the
identity 
$$
\nabla\cdot (h V\otimes V)=\nabla\cdot(hV)V+h(V\cdot\nabla)V,
$$
we thus get
$$
\begin{array}{r}
	\vspace{0.5em}
	(I+\mu\alpha h\cT\frac{1}{h})\dt (hV)
        + \eps(I+\mu \alpha h\cT\frac{1}{h})\nabla\cdot (hV\otimes V)\\
        + (I-\mu(1-\alpha)h\cT\frac{1}{h})h\nabla\zeta 
	+\eps\mu h\cQ_1(V)=0.
\end{array}
$$
The Green-Naghdi equations with improved dispersion can therefore be written in $(h,hV)$ variables
as
\be\label{eq6imphhV}
	\left\lbrace
	\begin{array}{l}
	\dsp \dt h +\eps \nabla\cdot (h V)=0,\\
       	\dsp (I+\mu\alpha h\cT\frac{1}{h})\dt (hV)
        +\eps(I+\mu \alpha h\cT\frac{1}{h})\nabla\cdot (hV\otimes V)\\
	\indent\indent \dsp +(I-\mu(1-\alpha)h\cT\frac{1}{h})h\nabla\zeta 
	+\eps\mu h\cQ_1(V)=0.
	\end{array}\right.
\ee
The second equation of (\ref{eq6imphhV}) requires the computation of third order derivatives of $\zeta$ that can be numerically stiff. It is however possible to show that these terms can be factorized by $I+\mu\alpha h\cT\frac{1}{h}$, up to a term involving only a first order derivative of $\zeta$:
$$
(I-\mu(1-\alpha)h\cT\frac{1}{h}) h \nabla\zeta = \frac{1}{\alpha} h \nabla \zeta + \frac{\alpha-1}{\alpha} (I+\mu\alpha h\cT\frac{1}{h}) h \nabla \zeta.
$$
The equations (\ref{eq6imphhV}) can therefore be reformulated as
\be\label{eq6imphhVbis}
	\left\lbrace
	\begin{array}{l}
	\dsp \dt h +\eps \nabla\cdot (h V)=0,\\
       	\dsp \dt (hV) 
        +\eps\nabla\cdot (hV\otimes V) + \frac{\alpha-1}{\alpha} h \nabla \zeta\\
	\indent\indent \dsp +(I+\mu\alpha h\cT\frac{1}{h})^{-1} [ \frac{1}{\alpha} h\nabla\zeta 
	+\eps\mu h\cQ_1(V) ] =0.
	\end{array}\right.
\ee
This formulation does not require the computation of any third-order derivative, allowing for more robust numerical computations, especially when the wave becomes steeper.

\subsection{Dimensionalized equations}\label{sectdim}
\label{eqdim}

Going back to variables with dimension, the system of equations (\ref{eq6imphhVbis}) reads
\be\label{GNdim}
\left\lbrace
\begin{array}{l}
\vspace{0.5em}
\dsp \dt h + \nabla\cdot (h V)=0,\\
\dsp \dt (hV) + \frac{\alpha-1}{\alpha} g h\nabla\zeta + \nabla\cdot (hV\otimes V)\\
\indent\indent\dsp + (I+\alpha h\cT\frac{1}{h})^{-1} [\frac{1}{\alpha}gh \nabla\zeta + h\cQ_1(V)] = 0,
\end{array}\right.
\ee
where the dimensionalized version of the operators $\cT$ and $\cQ_1$ correspond to (\ref{eq2}), (\ref{eq4}), (\ref{eq5}) and (\ref{eq11}) with $\beta=1$, and where $h$ now stands for the
water height with dimensions,
$$
h = h_0 + \zeta - b.
$$
Looking at the linearization of (\ref{GNdim}) around the rest state $h=h_0$, $V=0$, and flat bottom $b=0$, one derives
the dispersion relation associated to (\ref{GNdim}). It is found  by looking for plane wave solutions of the 
form $(\underline{h},\underline{h}\underline{V})e^{i({\bf k}\cdot x-\omega t)}$ to the linearized equations, and consists of
two branches parametrized by $\omega_{\alpha,\pm}(\cdot)$,
\begin{equation}\label{reldispGN}
\omega_{\alpha,\pm}({\bf k})=\pm \vert{\bf k}\vert \sqrt{gh_0}\sqrt{\frac{1+(\alpha-1)(\vert{\bf k}\vert h_0)^2/3}{1+\alpha(\vert{\bf k}\vert h_0)^2/3}}.
\end{equation}

\begin{figure}[t]
\centering
\psfrag{x}{}
\psfrag{y}[cr][bl][0.8]{$\dsp\frac{C_{GN}^{p/g}}{C_{S}^{p/g}}$}
\psfrag{Linear group velocity error}[cc][cl]{\hspace*{8em}\small \textit{Linear phase and group velocity errors}}
\hspace*{1em}\includegraphics[width=0.94\textwidth]{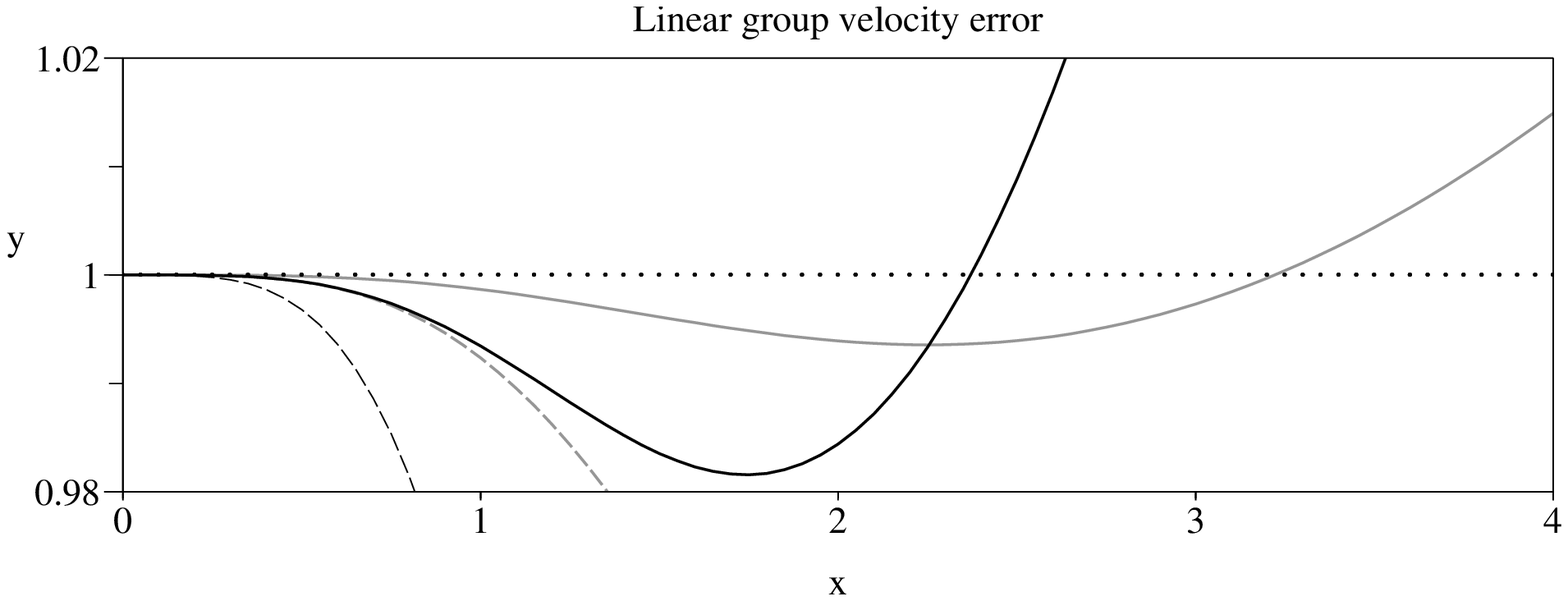}

\vspace*{1em}

\psfrag{x}[Bc][tc][0.85]{$kh_0$}
\psfrag{y}[cr][bl][1]{$\dsp\alpha_{opt}$}
\psfrag{Optimal value of a}[cc][cl]{\hspace*{5.5em}\small \textit{Local optimal values of $\alpha$}}
\hspace*{1em}\includegraphics[width=0.94\textwidth]{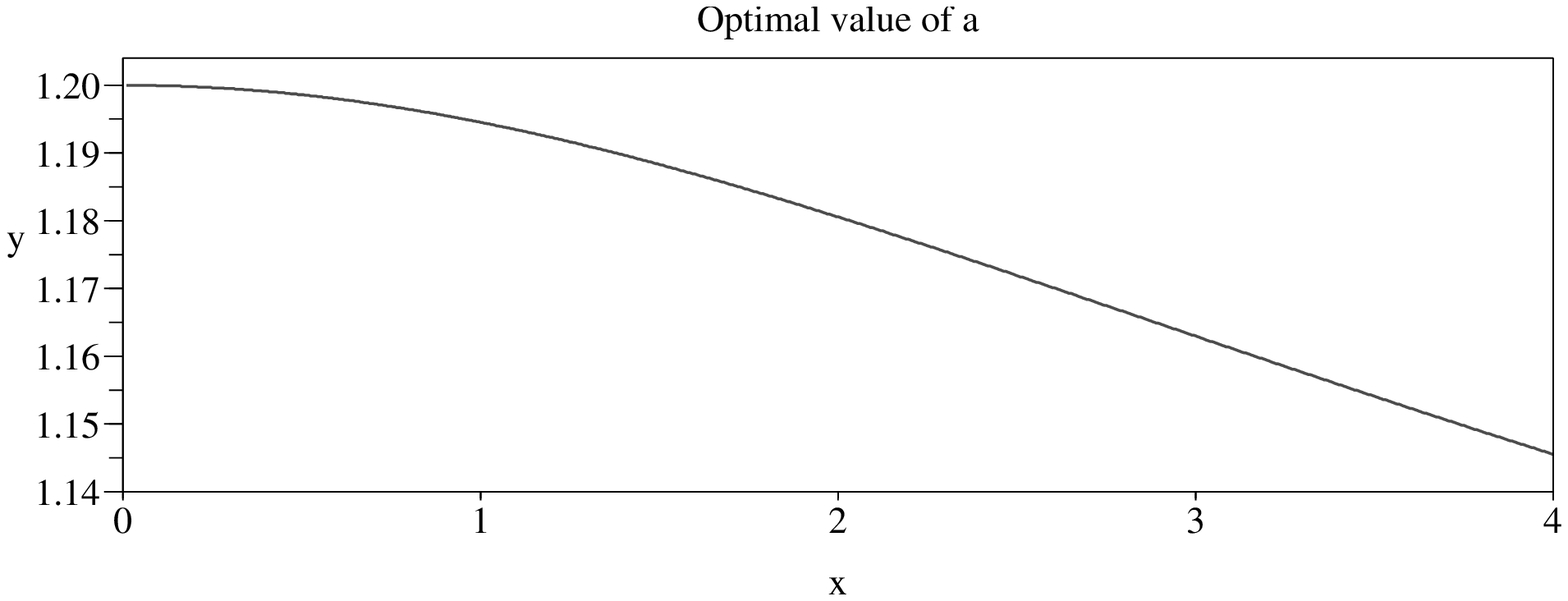}
\caption{Top: linear phase velocity (black) and group velocity (grey) errors for $\alpha=1.159$ (full line, global optimal value) and $\alpha=1$ (dashed line, original model). Bottom: local optimal values of $\alpha$ for $kh$ in $[0;4]$.}
\label{opt_alpha}
\end{figure}

As mentioned before, the parameter $\alpha$ can be helpful in optimizing the dispersive properties of the original model ($\alpha=1$), namely the linear phase and group velocities denoted by $C_{GN}^{\;\,p}$ and $C_{GN}^{\;\,g}$. By adjusting $\alpha$, we can minimize the error relative to the reference phase and group velocities $C_{S}^{p}$ and $C_{S}^{g}$ coming from Stokes~linear~theory. A classical approach consists in minimizing the averaged error over some range $kh_0 \in [0;K]$ (see \cite{DL}, \cite{Madsen} and \cite{cie2007} for further details). Here, minimizing the weighted\footnote{The squared relative error is weighted by $1/kh_0$ to keep the errors to a minimum for low wavenumbers.} averaged error over the range $kh_0 \in [0;4]$ yields the global optimal value $\alpha=1.159$, that is adopted, unless stated otherwise, throughout this paper\footnote{The dispersive correction used in \cite{cie2007} slightly differs from here, yielding a different definition of $\alpha$: denoting by $\widetilde{\alpha}$ the parameter used in \cite{cie2007}, the correspondence is given by $\widetilde{\alpha}=\frac{\alpha-1}{3}$.}. In Figure \ref{opt_alpha} (top), the ratios $C_{GN}^{\;\,p} / C_{S}^{p}$ and $C_{GN}^{\;\,g} / C_{S}^{p}$ are plotted against the relative water depth $kh_0$ for $\alpha = 1$ (original model) and $\alpha=1.159$ (optimized model).

\smallbreak

The global optimal value $\alpha=1.159$ is especially well-suited when considering irregular waves or regular waves over uneven bottoms, i.e. when multiple - or not easily predictable - wavelengths are involved. However, when considering monochromatic waves over flat bottoms, i.e. when only one wavelength is involved, $\alpha=1.159$ is not optimal anymore. In this particular case, an alternative approach consists in minimizing the error on the phase velocity for a specific value of $kh_0$, where $k$ corresponds to the wavenumber of the considered wave. For any discrete value of $kh_0$, one easily computes the corresponding optimal value of $\alpha$, denoted by $\alpha_{opt}(kh_0)$ and refered to as local optimal value.
In Figure \ref{opt_alpha} (bottom), $\alpha_{opt}$ is plotted against $kh_0$, for $kh_0 \in [0;4]$.

\section{Numerical methods}\label{sectNS}

We propose here to take advantage of the previous reformulation. First, this new formulation is well-suited for a splitting approach separating the hyperbolic and the dispersive part of
the equations (\ref{GNdim}). We present our splitting scheme in \S \ref{sectsplit}; we then
show in \S \ref{secthyp} and \S \ref{sectdisp} how we treat respectively the hyperbolic and dispersive parts
of the equations.

\subsection{The splitting scheme}
\label{sectsplit}

We decompose the solution operator $S(\cdot)$ associated to (\ref{GNdim}) at each time step
by the second order splitting scheme
\begin{equation}\label{S0}
	S(\delta_t)=S_1(\delta_t/2)S_2(\delta_t)S_1(\delta_t/2),
\end{equation}
where $S_1$ and $S_2$ are respectively associated to the hyperbolic and dispersive parts of
the Green-Naghdi equations (\ref{GNdim}). More precisely:
\begin{itemize}
\item  $S_1(t)$ is the solution operator
associated to NSWE
\be
	\label{S1}
	\left\lbrace
	\begin{array}{lcl}
	\vspace{0.25em}
	\dsp \dt h+\nabla\cdot (hV)&=&0,\\
	\dsp \dt (h V)+ \nabla (\frac{1}{2}g h^2) +\nabla\cdot (hV\otimes V) &=&-gh\nabla b.
	\end{array}\right.
\ee
\item $S_2(t)$ is the solution operator associated to the remaining (dispersive) part of the equations,
\be
	\label{S2}
	\left\lbrace
	\begin{array}{lcl}
	\dsp \dt h &=&0,\\
	\dsp \dt (h V) - \frac{1}{\alpha} g h \nabla \zeta + (I+\alpha h\cT\frac{1}{h})^{-1}\big[\frac{1}{\alpha} gh \nabla\zeta + h\cQ_1(V)\big]&=&0.
	\end{array}\right.
\ee
\end{itemize}

\begin{rmrk}
From this point, \emph{we only consider $1D$ surface waves}. The numerical implementation of our scheme for $2D$ surface
waves is left for future work.
\end{rmrk}

\begin{rmrk}
The friction term (see Remark \ref{friction}), when used, is included in $S_2(t)$.  
\end{rmrk}

Taking advantage of the hyperbolic structure of the NSWE, $S_1(t)$ is computed using a finite-volume approach, as described in the next subsection. As far as the operator $S_2(t)$ is concerned, we use a finite-difference approach, as shown in \S \ref{sectdisp}. 

\medbreak

Such a mixed finite-volume finite-difference method implies to work both on cell-averaged and nodal values for each unknown. 
We use the following notations:
\begin{nota}
- The numerical one-dimensional domain $\Omega$ is uniformly divided into $N$ cells $(C_i)_{1 \le i \le N}$ such that $C_i = [x_{i-1}, x_i]$, where  $(x_i)_{0 \le i \le N}$ are the $N+1$ nodes of the regular grid.\\
-  We denote by $\delta_x$ the cell size and by $\delta_t$ the time step.\\
- We write $w_i^n$  the nodal value of $w$ at the $i^{th}$ node $x_i$ and at time $t_n = n\delta_t$.\\
- We denote by $\bar{w}_i^n$ the averaged value of $w$ on the $i^{th}$ cell $C_i$ at time $t_n = n\delta_t$.
\end{nota}

\medbreak

The choice of a finite difference method for solving $S_2$ also entails to switch between the cell-averaged values and the nodal values of each unknown, in a suitable way that preserves the global spatial order of the scheme. Using classical fourth-order Taylor expansions, we easily recover the following relations that allow to switch between the finite volume unknowns $(\bar{w}_i^n)_{1 \le i \le N}$ and the finite difference unknowns $(w_i^n)_{0 \le i \le N}$ at each time step:
\begin{equation}\label{P1}
\frac{1}{6}w_{i-1} + \frac{2}{3} w_i + \frac{1}{6}w_{i+1} = \frac{1}{2}(\bar{u}_i + \bar{u}_{i+1}) + O(\delta_x^4),\;0 \le i \le N,
\end{equation}
and
\begin{equation}\label{P2}
\bar{u}_i = -\frac{1}{24}w_{i-2} + \frac{13}{24}w_{i-1} + \frac{13}{24}w_i - \frac{1}{24}w_{i+1} + O(\delta_x^4),\;1 \le i \le N,
\end{equation}
with adaptations at the boundaries following the method presented in \S \ref{bound}.\\
We can easily check that (\ref{P1}), (\ref{P2}) preserves the steady state at rest, and that these formulae are precise up to $O(\delta_x^4)$ terms, thus preserving the global order of the scheme.

\subsection{Spacial discretization of the hyperbolic component $S_1(\cdot)$}
\label{secthyp}

 When specified in one space dimension, the
system under consideration reads as follows:
\be\label{stvenant}
\left\{
\begin{array}{l}
\dt h + \dx(hu) = 0,\\
\dt(hu) + \dx ( hu^2 + gh^2/2 ) = -gh\partial_x b
\end{array}\right.
\ee
For the sake of simplicity in the notations, it is convenient to
 rewrite the system (\ref{stvenant}) in the following condensed form:
\be\label{stvenant2d2}
\dt\vw + \dx\vf(\vw) = \vS(\vw,b),
\ee
with
\be
\vw = \left(\begin{array}{c}
h\\
hu
\end{array}\right),
\quad
\vf(\vw) =\left( \begin{array}{c}
hu\\
hu^2 + \dsp\fr g 2 h^2
\end{array}\right)
\quad\mbox{ and }\quad
\vS(\vw) =\left( \begin{array}{c}
0\\
-gh\partial_x b
\end{array}\right),
\ee
where $\vw: \R\,\times\,\R^+ \rightarrow \Omega$ is the state
vector in conservative variables and $\vf(\vw):  \Omega
\rightarrow \R^2$ stands for the flux function. The convex set $\Omega$ of
the admissible states is defined by
$$
\Omega = \left\{ \vw\,\in\,\R^2;\;h\,\geq\,0, u\,\in\,\R \right\}.
$$
Considering numerical approximations of system (\ref{stvenant2d2}), we seek a numerical scheme that provides stable simulations of the processes occurring in surf and swash areas, with a precise control of the spurious effects induced by numerical dissipation and dispersion. Moreover, the scheme should be able to handle the complex interactions between waves and topography, including the preservation of motionless steady states:
$$
u=0,\quad h+b = {\rm cste}.
$$
In this way, we use a low-dissipation and well-balanced extension of the robust finite volume scheme introduced in \cite{berthon_marche}. The main features of the first order scheme are recalled in \S \ref{FO}, its higher-order and well-balanced extension presented respectively in \S \ref{HE} and \S \ref{WB}.

\subsubsection{First order finite-volume scheme for the homogeneous system}\label{FO}

The homogeneous NSWE associated with (\ref{stvenant}), given by
\bea\label{stvh}
\dt\vw+\dx\vf(\vw)=0,
\eea
is known to be hyperbolic over $\Omega$. As a consequence, the
solutions may develop shock discontinuities. In order to rule out the
unphysical solutions, the system (\ref{stvh}) must be supplemented by
entropy inequalities (see for instance \cite{bouchut} and references therein). \\
The spatial discretization of the homogeneous system (\ref{stvh})  can be recast under the following classical semi-discrete finite-volume formalism:
$$
\frac{d}{dt}\bar{\vw}_i (t) +   \frac{1}{\delta_x}   \Bigl( \tilde{\vf}  \bigl(\bar{\vw}_{i},  \bar{\vw}_{i+1}\bigr)   -  \tilde{\vf} \left(\bar{\vw}_{i-1},\bar{\vw}_{i} \right) \Bigr)= 0
$$
where $\tilde{\vf}$ is a numerical flux function based on a conservative flux consistent with the homogeneous NSWE. For the numerical validations shown in \S\ref{numval}, we use the  numerical flux issued from the relaxation approach introduced in \cite{berthon_marche}.

\begin{rmrk}
The robustness of this finite volume scheme for the homogeneous NSWE is shown in \cite{berthon_marche}, where the detailed study of the relaxation approach is performed.
\end{rmrk}

\subsubsection{A robust high-order extension}\label{HE}

To reduce both numerical dissipation and dispersion within the hyperbolic component $S_1(\cdot)$, high order reconstructed states at each
interface have to be considered. Following the classical MUSCL approach \cite{vanleer}, we consider the modified scheme:
\be\label{scheme2}
\frac{d}{dt}\bar{\vw}_i (t) + \frac{1}{\delta_x}\left( \tilde{\vf} \left( \bar{\vw}_{i,r}^n,\bar{\vw}_{i+1,l}^n
                    \right)- \tilde{\vf} \left(\bar{\vw}_{i-1,r}^n,\bar{\vw}_{i,l}^n
                    \right) \right)=0,
\ee
where $\bar{\vw}_{i,l}^n$ and $\bar{\vw}_{i,r}^n$ are high-order interpolated values of the cell-averaged solution, respectively at the left and right interfaces of the cell $C_i$.  The low dissipation reconstruction proposed in \cite{camarri} is used. Considering a cell $C_i$,  and the corresponding constant value $\bar{h}_i^n$, we introduce linear reconstructed left and right values $\bar{h}_{i,l}^n$ and  $\bar{h}_{i,r}^n$ as follows:

\beq\label{interp}
\bar{h}_{i,r}^n = \bar{h}_i^n + \frac{1}{2} \overline{\delta h}_{i,r}^n \mbox{\;\;\; and \;\;\;} \bar{h}_{i,l}^n = \bar{h}_i^n - \frac{1}{2}\overline{\delta h}_{i,l}^n.
\eeq
The corresponding gradients are built following the five points stencil:
\beq\label{interp_plus}
\begin{split}
\overline{\delta h}_{i,r}^n = (1- & \nu)(\bar{h}_{i+1}^n - \bar{h}_i^n) +  \nu (\bar{h}_{i}^n - \bar{h}_{i-1}^n) \\ &+ \xi^c(-\bar{h}_{i-1}^n + 3 \bar{h}_i^n - 3 \bar{h}_{i+1}^n +\bar{h}_{i+2}^n) \\&+ \xi^d(-\bar{h}_{i-2}^n + 3 \bar{h}_{i-1}^n -3 \bar{h}_{i}^n + \bar{h}_{i+1}^n),
\end{split}
\eeq
\beq
\begin{split}\label{interp_min}
\overline{\delta h}_{i,l} = (1- & \nu)(\bar{h}_{i}^n - \bar{h}_{i-1}^n) +  \nu (\bar{h}_{i+1}^n - \bar{h}_{i}^n) \\ &+ \xi^c (-\bar{h}_{i-2}^n + 3 \bar{h}_{i-1}^n -3 \bar{h}_{i}^n + \bar{h}_{i+1}^n)\\&+ \xi^d(-\bar{h}_{i-1}^n + 3 \bar{h}_i^n - 3 \bar{h}_{i+1}^n +\bar{h}_{i+2}^n),
\end{split}
\eeq
and the coefficients $\nu$, $\xi^c$ and $\xi^d$ are set respectively to $\fr 1 3$, $-\frac{1}{10}$ and $-\frac{1}{15}$, leading 
to better dissipation and dispersion properties in the truncature error.

\medbreak

When the generation of shock waves occurs during computation, the previous interpolation has to be embedded into a limitation procedure to keep the scheme non oscillatory and positive. We suggest to use a three-entry limitation, especially designed to generate a positive scheme of higher possible order far from extrema and discontinuities. Scheme (\ref{scheme2}) thus becomes
\be\label{scheme3}
\frac{d}{dt}\bar{\vw}_i (t) + \frac{1}{\delta_x}\left( \tilde{\vf}\left({}^L\bar{\vw}_{i,r}^n,{}^L\bar{\vw}_{i+1,l}^n
                    \right)-\tilde{\vf} \left( {}^L\bar{\vw}_{i-1,r}^n,{}^L\bar{\vw}_{i,l}^n
                    \right) \right) = 0.
\ee
The limited high-order reconstructed values are defined, considering for instance the water height $h$, as
\beq
{}^L\bar{h}_{i,r}^n = \bar{h}_i^n + \frac{1}{2}L_{i,r}(\bar{h}^n) \mbox{\;\;\; and \;\;\;} {}^L\bar{h}_{i,l}^n = \bar{h}_i^n - \frac{1}{2}L _{i,l}(\bar{h}^n).
\eeq
To define $L_{i,r}(\bar{h}^n)$ and $L _{i,l}(\bar{h}^n)$, we use the following limiter:
\beq\label{superbee}
 L(u,v,w) = \left\{ \begin{array}{ll}
    0 & \mbox{ if } uv\leq 0,\\
    {\it \mbox{sign}}(u)\,\min(2\vert u \vert, 2\vert v \vert,w) & \mbox{otherwise}.
   \end{array}\right.\\
\eeq 
Relying on (\ref{superbee}), we then define the limiting process as
$$
 L_{i,r}(\bar{h}^n)= L(\overline{\delta h}_i^{n,-},   \overline{\delta h}_i^{n,+},  \overline{\delta h}_{i,r}^n) \mbox{\;\; and\;\;}   L _{i,l}(\bar{h}^n)= L(\overline{\delta h}_i^{n,+}, \overline{\delta h}_i^{n,-},  \overline{\delta h}_{i,l}^n),
$$
where $\overline{\delta h}_i^{n,+} = \bar{h}_{i+1}^n - \bar{h}_i^n$ and $\overline{\delta h}_i^{n,-} = \bar{h}_{i}^n - \bar{h}_{i-1}^n$ are upstream and downstream variations, and $\overline{\delta h}_{i,r}^n$ and $\overline{\delta h}_{i,l}^n$ taken from (\ref{interp_plus}) and (\ref{interp_min}).\\
Such limited high order reconstructions must also be performed for the other conservative variable $hu$.

\begin{rmrk}
It is straightforward that when the considered conservative variable is smooth enough, this limiter preserves the high order accuracy of the reconstructions (\ref{interp_plus}) and (\ref{interp_min}). In addition, this high order reconstruction and the limitation process can easily be extended to non-uniform meshes. 
\end{rmrk}

\begin{rmrk}
The robustness of the resulting high order relaxation scheme can be proved following the lines of \cite{berthon_marche}.
\end{rmrk}

\subsubsection{Well-balancing for steady states}\label{WB}
We finally introduce a well-balanced discretization of the topography source term. Scheme (\ref{scheme3}) is embedded within a {\it hydrostatic reconstruction} step \cite{bouchut}.

To achieve both well-balancing and high order accuracy requirements, we have to consider not only high order reconstructions of the conservative variables, as done in \S \ref{HE}, but also of the surface elevation $\zeta$. The resulting finite volume scheme is able to preserve both motionless steady states and water height positivity. The reader is referred to  \cite{bouchut} for a detailed study of the {\it hydrostatic reconstruction} method, including robustness, stability and semi-discrete entropy inequality results.

\subsection{Spacial discretization of the dispersive component  $S_2(\cdot)$}
\label{sectdisp}

The system corresponding to the operator $S_2(\cdot)$ writes in one dimension
\be
	\label{S2_1D}
	\left\lbrace
	\begin{array}{lcl}
	\dsp \dt h &=&0,\\
	\dsp \dt (h u) - \frac{1}{\alpha} g h \partial_x \zeta + (1+\alpha h\cT\frac{1}{h})^{-1}\big[\frac{1}{\alpha} gh \partial_x\zeta +
        h\cQ_1(u)\big]&=&0
	\end{array}\right.
\ee
where the operators $\cT$ and $\cQ_1$ are explicitly given by
\be
\cT w = -\frac{h^2}{3} \partial_x^2 w - h \partial_x h \partial_x w + (\partial_x \zeta \partial_x b + \frac{h}{2} \partial_x^2 b) w,
\ee
and
\be
\cQ_1(u) = 2h\partial_x (h+\frac{b}{2}) (\partial_x u)^2 + \frac{4}{3} h^2 \partial_x u \partial_x^2 u + h \partial_x^2 b u \partial_x u
+ (\partial_x \zeta \partial_x^2 b + \frac{h}{2} \partial_x^3 b) u^2.
\ee
As specified in \S \ref{sectsplit}, the system (\ref{S2_1D}) is solved at each time step using a classical finite-difference technique. The spatial derivatives are discretized using the following fourth-order formulae:
\begin{eqnarray}
(\delta_x w)_i & = & \frac{1}{12 \delta_x} (-w_{i+2} + 8 w_{i+1} - 8 w_{i-1} + w_{i-2}),\nonumber\\
(\partial_x^2 w)_i & = & \frac{1}{12 \delta_x^2} (-w_{i+2} + 16 w_{i+1} - 30 w_i + 16 w_{i-1} - w_{i-2}),\nonumber\\
(\partial_x^3 w)_i & = & \frac{1}{8 \delta_x^3} (-w_{i+3} + 8 w_{i+2} - 13 w_{i+1} + 13 w_{i-1} - 8 w_{i-2} + w_{i-3}).\nonumber
\end{eqnarray}
Boundary conditions are imposed using the method presented in \S \ref{bound}.

\subsection{Time discretization and dispersive properties}\label{secttime}

\subsubsection{Time discretization}
As far as time discretization is concerned, we choose to use explicit methods. The systems corresponding to $S_1$ and $S_2$ are integrated in time using a classical fourth-order Runge-Kutta scheme.

\subsubsection{Dispersive properties}\label{reldispdt}

We now turn to investigate the dispersive properties of our numerical scheme. Since the main originality of this approach is the 
splitting in time of the hyperbolic and dispersive parts, we consider here the semi-discretized in time version of our
numerical scheme. An extension to the fully discretized scheme is of course possible, but extremely technical, and
would not bring any significant 
insight on the dispersive properties of the hyperbolic/dispersive splitting.

We recall that the dispersion relation associated to the Green-Naghdi (with improved frequency dispersion) equations 
(\ref{GNdim}) is given by (\ref{reldispGN}) or, for $1D$ surface waves,
$$
\omega_{\alpha,\pm}(k)=\pm{k}\sqrt{gh_0}\sqrt{ \frac{1+(\alpha-1)(kh_0)^2/3}{1+\alpha(k h_0)^2/3}}.
$$
The dispersion relation corresponding to our semi-discretized (in time) splitting scheme is given by the following
proposition.
\begin{proposition}\label{propdisp}
The dispersion relation associated to the semi-discretized \\scheme (\ref{S0}), (\ref{S1}), (\ref{S2}) is given by
$$
\omega_{sd,\pm}(k)=\omega_{\alpha,\pm}(k)+\frac{\delta_t^2}{24}\omega_{\alpha,\pm}(k)^3\Big(\frac{(kh_0)^2}{3+(\alpha -1)(kh_0)^2}\Big)^2+O(\delta_t^3).
$$
\end{proposition}
\begin{rmrk}
The proposition above shows that the semi-discretized dispersion relation approaches the exact dispersion relation of the
Green-Naghdi equations (\ref{GNdim}) at order $2$ in $\delta_t$. An additional information is that the $O(\delta_t^2)$
error made by the splitting scheme is always real. Therefore, the numerical errors are of \emph{dispersive} type and there
is no linear instability induced by the splitting.
\end{rmrk}
\begin{rmrk}
Since the main error in the dispersive relation is of dispersive type, it is natural
to try to remove it with techniques inspired by the classical Lax-Wendroff scheme.
This is possible, but this does not yield better results than the $\delta_t$-optimized
choice of the frequency parameter $\alpha$ (see below), which is a much simpler method. 
\end{rmrk}
\begin{proof}
For the sake of clarity, we still denote by $S_1(\cdot)$ and $S_2(\cdot)$ the solution operators associated to the 
semi-discretized version of the linearization of (\ref{S1}) and (\ref{S2}) around the rest state (and flat bottoms).\\
{\bf Step 1.} We show here that
$$
\forall {\bf w}^0\in \R^2, \qquad S_1(\delta_t)({\bf w}^0 e^{ikx})=
\left(\begin{array}{cc} \alpha_1(\delta_t) & \alpha_2(\delta_t)ik \\
gh_0 \alpha_2(\delta_t)ik & \alpha_1(\delta_t)\end{array}\right){\bf w}^0 e^{ikx},
$$
with 
$$
\alpha_1(\delta_t)=1+\frac{\delta_t^2}{2}(-gh_0 k^2)+\frac{\delta_t^4}{24}(-gh_0k^2)^2,\qquad
\alpha_2(\delta_t)=\delta_t +\frac{\delta_t^3}{6}(-gh_0k^2).
$$
Since the linearization of (\ref{S1}) around the rest state and flat bottom can be written in compact form as
$$
\dt {\bf w}+A(\partial_x){\bf w}=0,\quad \mbox{ with }\quad
A(\partial_x)=\left(\begin{array}{cc} 0 & \partial_x\\ gh_0\partial_x & 0\end{array}\right),
$$
the quantity $S_1(\delta_t)({\bf w}^0 e^{ikx})$ corresponding to the RK4 time discretization is given by
$$
S_1(\delta_t)({\bf w}^0 e^{ikx})=\Big(1+\delta_t A(ik)+\frac{\delta_t^2}{2}A(ik)^2+\frac{\delta_t^3}{6}A(ik)^3+\frac{\delta_t^4}{24}A(ik)^4\Big){\bf w}^0 e^{ikx}. 
$$
A simple computation thus yields the result.\\
{\bf Step 2.} We show here that
$$
\forall {\bf w}^0\in \R^2, \qquad S_2(\delta_t)({\bf w}^0 e^{ikx})=
\left(\begin{array}{cc} 1 & 0\\ gh_0\gamma ik \delta_t & 1 \end{array}\right){\bf w}^0 e^{ikx},
$$
with
$$
\gamma=-\frac{(kh_0)^2/3}{1+\alpha(kh_0)^2/3}.
$$
Since the linearization of (\ref{S2}) around the rest state and flat bottom can be written in compact form as 
$$
\dt {\bf w}+B(\partial_x){\bf w}=0,\quad\mbox{ with }\;
B(\partial_x)=\left(\begin{array}{cc} 0 & 0\\ -gh_0(1-\frac{\alpha}{3}h_0^2\partial_x^2)^{-1} (-\frac{1}{3}h_0^2\partial_x^2)\partial_x& 0\end{array}\right),
$$
the quantity $S_2(\delta_t)({\bf w}^0 e^{ikx})$ corresponding to the RK4 time discretization is given by
$$
S_2(\delta_t)({\bf w}^0 e^{ikx})=\big(1+\delta_t B(ik)\big){\bf w}^0 e^{ikx},
$$
where we used the fact that $B(ik)^2=0$. The result follows directly.\\
{\bf Step 3.} By a direct computation, we get that
$$
\forall {\bf w}^0\in \R^2, \qquad S_1(\delta_t/2)S_2(\delta_t)S_1(\delta_t/2)({\bf w}^0 e^{ikx})=
(I+\delta_t \,M ){\bf w}^0 e^{ikx},
$$
with $M=(m_{ij})_{1\leq i,j\leq 2}$ given by
\begin{eqnarray*}
m_{11}=m_{22}&=& -\frac{gh_0}{2}(1+\gamma)k^2\delta_t+O(\delta_t^3),\\ 
m_{12}&=&
ik-i\frac{gh_0}{6}(1+\frac{3}{2}\gamma)k^3\delta_t^2+O(\delta_t^4),\\
m_{21}&=& i gh_0(1+\gamma)k-i\frac{(gh_0)^2}{6}(1+\frac{3}{2}\gamma)k^3\delta_t^2+O(\delta_t^4) .
\end{eqnarray*}
{\bf Step 4.} End of the proof. We deduce from the
previous steps that if ${\bf w}^0e^{ikx-\omega t}$ is a plane wave solution for the semi-discretized scheme, then one has
$$
e^{-i\omega \delta_t}{\bf w}^0=(I+\delta_t\, M) {\bf w}^0,
$$
and $\frac{e^{-i\omega \delta_t}-1}{\delta_t}$ is therefore an eigenvalue of $M$. After some simple computations, we thus
get
\begin{equation}\label{reld}
\frac{e^{-i\omega \delta_t}-1}{\delta_t}=\lambda_\pm,
\end{equation}
with
\begin{equation}\label{reld1}
\lambda_\pm=-i\omega_{\alpha,\pm}(k)-\frac{1}{2}\omega_{\alpha,+}^2(k)\delta_t+\frac{i}{24} \omega_{\alpha,\pm}(k)^3(4-\frac{\gamma^2}{(1+\gamma)^2})\delta_t^2+O(\delta_t^3).
\end{equation}
By identifying the Taylor expansion of the left-hand-side of (\ref{reld}) with (\ref{reld1}), we deduce that
\begin{equation}\label{reld2}
\omega=\omega_{\alpha,\pm}+\frac{1}{24}\omega_{\alpha,\pm}(k)^3\frac{\gamma^2}{(1+\gamma)^2}\delta_t^2+O(\delta_t^3),
\end{equation}
and the result follows.
\end{proof}

Starting from the previous expression and dropping the $O(\delta_t^3)$ term, one easily obtains the semi-discrete linear phase and group velocities $C_{GN}^{\;\,p}(\delta_t)$ and $C_{GN}^{\;\,g}(\delta_t)$, and computes the semi-discrete error relative to the reference velocities $C_{S}^{p}$ and $C_{S}^{g}$. Obviously, the global value $\alpha=1.159$ is no longer optimal with the additional $O(\delta_t^2)$ term, and we need to compute new optimal values of $\alpha$ that depend on $\delta_t$. As in \S \ref{sectdim} and for discrete values of the non-dimensional time step $\widetilde{\delta_t} := \sqrt{\frac{g}{h0}}\delta_t$, we look for 1) the global optimal value $\overline{\alpha_{opt}}(\delta_t)$ over the range $[0;3]$, and 2) the local optimal values $\alpha_{opt}(\delta_t,kh_0)$ for some discrete values~of~$kh_0$.

\vspace{1em}

\begin{figure}[t]
\centering
\psfrag{x}[Bc][tc]{}
\psfrag{y}[cr][bl][1]{$\dsp\overline{\alpha_{opt}}$}
\psfrag{Optimal value of a}[cc][cl]{\hspace*{5.5em}\small\textit{Global optimal values of $\alpha$}}
\hspace*{1em}\includegraphics[width=0.97\textwidth]{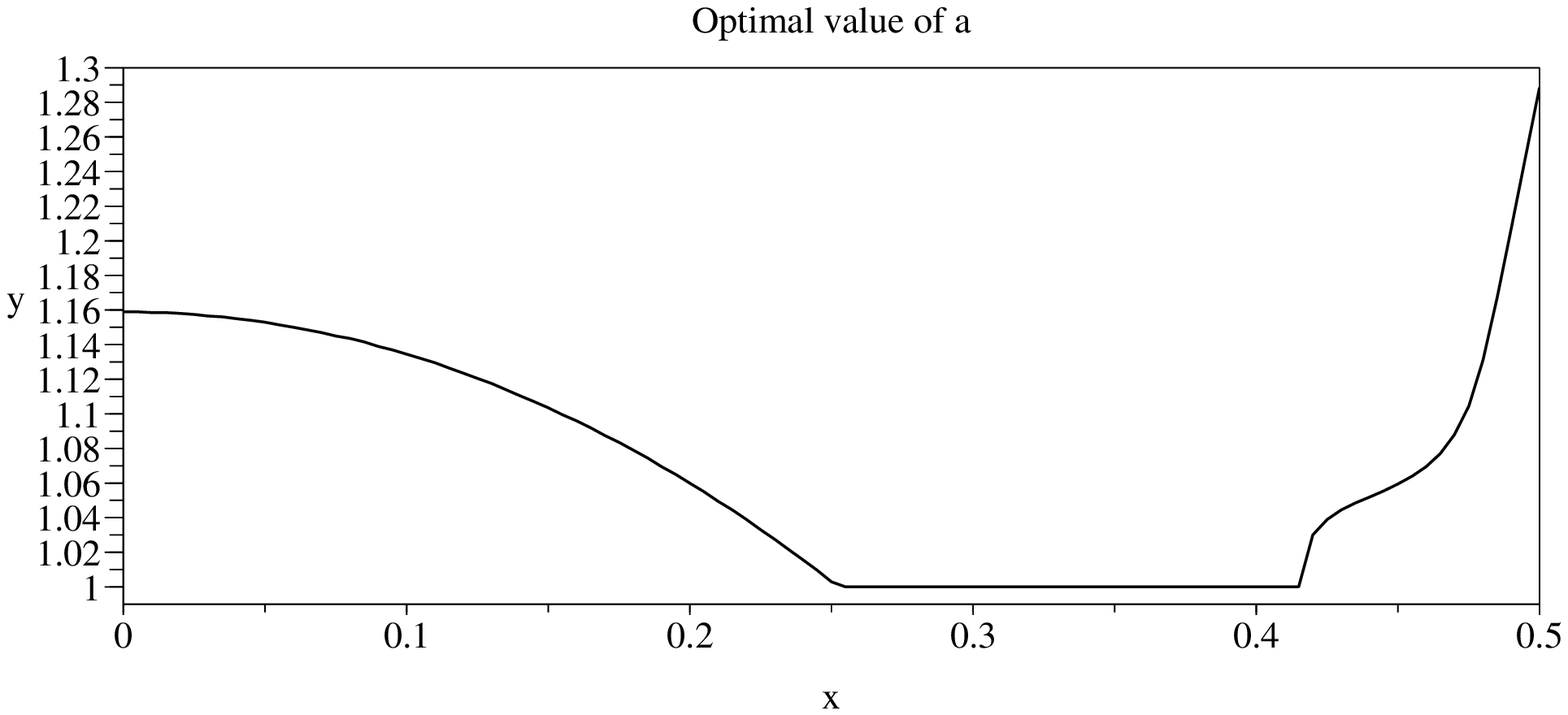}

\vspace*{0.75em}

\psfrag{x}[Bc][tc][0.85]{$\sqrt{g/h_0}\,\delta_t$}
\psfrag{y}[cr][bl][1]{$\dsp\alpha_{opt}$}
\psfrag{Optimal value of a}[cc][cl]{\hspace*{5.7em}\small \textit{Local optimal values of $\alpha$}}
\hspace*{1em}\includegraphics[width=0.97\textwidth]{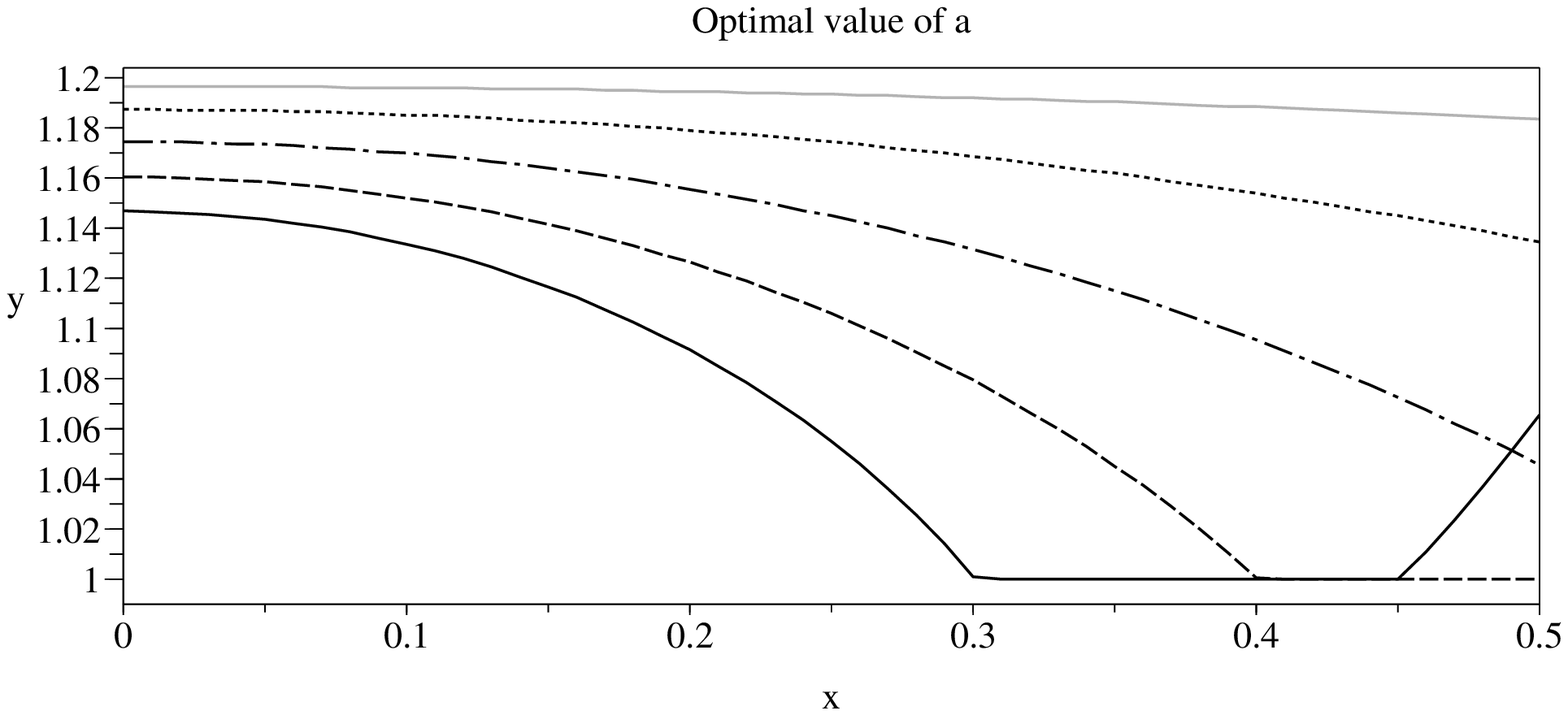}

\vspace{0.75em}

\caption{Top: global optimal values of $\alpha$ against $\widetilde{\delta_t}$. Bottom: local optimal values of $\alpha$ against $\widetilde{\delta_t}$ for $kh_0=\pi/4$ (full grey line), $kh_0=\pi/2$ (dotted line), $kh_0=3\pi/4$ (dash-dotted line), $kh_0=\pi$ (dashed line), $kh_0=5\pi/4$ (full black line).}
\label{opt_alpha_dt}
\end{figure}

Results are gathered in Figure \ref{opt_alpha_dt}: the top figure plots the global optimal value $\overline{\alpha_{opt}}$ against $\widetilde{\delta_t}$, while the bottom figure plots the local optimal values $\alpha_{opt}$ against $\widetilde{\delta_t}$, each curve corresponding to a discrete value of $kh_0$. We point out that in both approaches, the computed optimal value of $\alpha$ was sometimes found lower than $1$, for instance when $\widetilde{\delta_t} \in [0.26;0.42]$ for the global optimal value $\overline{\alpha_{opt}}$. Since taking $\alpha < 1$ induces some high-frequency instabilities, the optimal value of $\alpha$ has been taken equal to $1$ in such cases. However, it is worth remarking that for these problematic $\widetilde{\delta_t}$-regions, the model that provides the best dispersive properties - among the stable ones -  is the original Green-Naghdi model. 

We finally refer to \S \ref{fenton} for numerical simulations showing  the consequences of our choice to 
optimize $\alpha$ taking into account the dispersive effects of the time discretization.

\subsection{Boundary conditions}\label{bound}

The boundary conditions for the hyperbolic part $S_1$ of the splitting are treated as in \cite{marche_bonneton}. More precisely, as the simulations shown in this work do not require complex Riemann invariants based inflow, outflow or absorbing conditions, we simply introduce "ghosts cells" respectively at left and right boundaries of the domain, and suitable relations are imposed on the cell-averaged quantities~:
\begin{itemize}
\item $\bar{w}_{-k+1} = \bar{w}_{N-k+1}$ and $\bar{w}_{N+k} = \bar{w}_{k}$, $k \ge 1$, for periodic conditions on the left and right boundaries,
\item $\bar{w}_{-k+1} = \bar{w}_{k}$ and $\bar{w}_{N+k} = \bar{w}_{N-k+1}$, $k \ge 1$, for homogeneous Neumann conditions on the left and right boundaries,
\item $\bar{w}_{-k+1} = -\bar{w}_{k}$ and $\bar{w}_{N+k} = -\bar{w}_{N-k+1}$, $k \ge 1$, for homogeneous Dirichlet conditions on the left and right boundaries.
\end{itemize}
For the dispersive part $S_2$ of the splitting, the boundary conditions are simply imposed by reflecting - periodically for periodic conditions, evenly for Neumann conditions and oddly for Dirichlet conditions - the coefficients associated to stencil points that are located outside of the domain. The advantage of this method is to avoid introducing decentered formulae at the boundaries, while maintaining a regular structure in the discretized model:
\begin{itemize}
\item $w_{-k} = w_{N-k}$ and $w_{N+k-1} = w_{k-1}$, $k \ge 1$, for periodic conditions on the left and right boundaries,
\item $w_{-k} = w_{k}$ and $w_{N+k} = w_{N-k}$, $k \ge 1$, for homogeneous Neumann conditions on the left and right boundaries,
\item $w_{-k} = -w_{k}$ and $w_{N+k} = -w_{N-k}$, $k \ge 1$, for homogeneous Dirichlet conditions on the left and right boundaries.
\end{itemize}
Solid wall effects on the left or right boundary can be easily reproduced by imposing an homogeneous Neumann condition on $h$ and $\bar{h}$, and an homogeneous Dirichlet condition on $hu$ and $\bar{hu}$, with the previous methods.\\
When the water depth vanishes, a small routine is applied to ensure stability on the results: on each cell, if $h$ is smaller than some threshold $\epsilon$ then we impose the values $h=\epsilon$ and $v=0$.

\subsection{Wave breaking}
In order to handle wave breaking, we switch from the Green-Naghdi equations to the NSWE, locally in time and space, by skipping the dispersive step $S_2(\delta t)$ when the wave is ready to break. In this way, we only solve the hyperbolic part of the equations for the wave fronts, and the breaking wave dissipation is represented by shock energy dissipation (see also  \cite{bonn2007}, \cite{marche_bonneton} and \cite{brocchini2}). \\ 
To determine where to suppress the dispersive step at each time step, we use the first half-time step $S_1$ of the time-splitting as a predictor to assess the local energy dissipation, given by
\begin{equation}
 {D_i} =  -(\partial_t {\cal E} + \partial_x {\cal F}),
\end{equation}
with $\cal{E}$ = $\frac{\rho}2 (h u^2 + g \zeta^2)$ the energy density and
$ \cal{F}$= $\rho h u (\frac{u^2}2 + g\zeta )$ the energy flux density.
 This dissipation is close to zero in regular wave regions, and forms a peak when shocks are appearing. We can then easily locate the eventual breaking wave fronts at each time step, and skip the dispersive step only at the wave fronts.

\section{Numerical validation}\label{numval}

\subsection{Propagation of a solitary wave}\label{solwav}

It is known that for horizontal bottoms, the Green-Naghdi model with $\alpha=1$ have an exact solitary wave solution given by
\begin{equation}
\left\{
\begin{array}{l}
\vspace{0.2em}
\dsp h(x,t)=h_0+a\hspace{0.05cm}\mbox{sech}^2(\kappa(x-ct)),\\
\vspace{0.2em}
\dsp u(x,t)=c\Big{(}1-\frac{h_0}{h(x,t)}\Big{)},\\
\vspace{0.3em}
\dsp \kappa=\frac{\sqrt{3H}}{2h_0\sqrt{h_0+H}},
\qquad c=\sqrt{g(h+H)},
\end{array}
\right.
\label{solw}
\end{equation}
This family of solutions can be used as a validation tool for our present numerical scheme. We successively consider the propagation of two solitary waves of different relative amplitude $a/h_0$, on a $30\,m$ long domain with a constant depth $h_0 = 0.5\,m$. The considered relative amplitudes are $a/h_0 = 0.05$ for the first solitary wave, and $a/h_0 = 0.2$ for the second one. Periodic conditions are imposed on each boundary, and the initial surface and velocity profiles are centered at the middle of the domain.

\medbreak

\noindent In order to assess the convergence of our numerical scheme, the numerical solution is computed for several time steps $\delta_t$ and cell sizes $\delta_x$, over a sufficient duration $T=3\,s$. Starting with $\delta_x = 1\,m$ and $\delta_t = \delta_x / \sqrt{gh_0} = 0.45\,s$, we successively divide the time step by two, while keeping the CFL equal to $1$. For each computation and each discrete time $t_n = n \delta_t$, the relative errors $E_\zeta^n$ and $E_u^n$ on the free surface elevation and the averaged velocity are computed using the discrete $L^{\infty}$ norm 
$||.||_{\infty}$ :
$$
E_\zeta^n = \frac{||h_{num}-h_{sol}||_{\infty}}{||h_{sol}-h_0||_{\infty}}\quad;\quad
E_u^n = \frac{||u_{num}-u_{sol}||_{\infty}}{||u_{sol}||_{\infty}}
$$
where $(h_{num},u_{num})$ are the numerical solutions and $(\zeta_{sol},u_{sol})$ denotes the analytical ones coming from (\ref{solw}).

\medbreak

\begin{figure}[!t]
\psfrag{dt}[c][c]{$\delta_t \;\mbox{or}\; \delta_x / \sqrt{gh_0}$}
\psfrag{err}[cc][cl][1][90]{$\max{E_\zeta}$}
\includegraphics[width=0.95\textwidth]{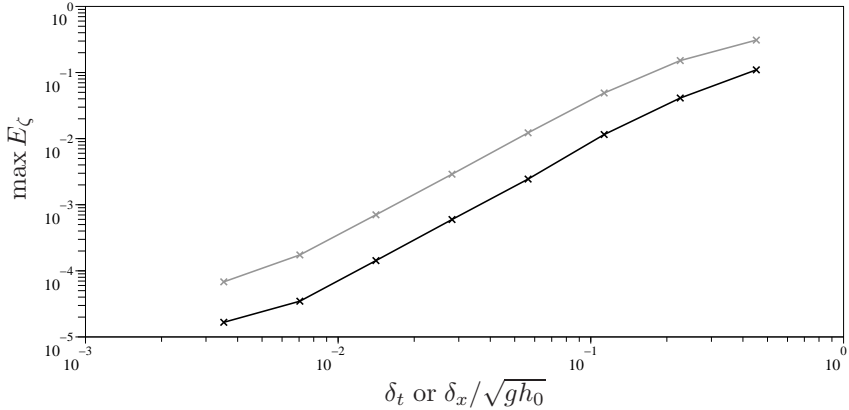}
\caption{Propagation of a solitary wave over a flat bottom, maximum of the relative error on the free surface elevation. Case $a/h_0 = 0.05$ in black, case $a/h_0 = 0.2$ in grey.}
\label{valid}
\end{figure}

Results are gathered in Figure \ref{valid}, where $\max{E_\zeta}$ is plotted against $\delta_t$, for the two considered relative amplitudes $a/h_0 = 0.05$ and $a/h_0 = 0.2$. In both cases, the convergence of our numerical scheme is clearly demonstrated. Furthermore, computing a linear regression on all points yields a slope equal to $1.91$ for the first case and $1.83$ for the second one. This result is coherent since the global order of our scheme is obviously limited by the order of the splitting method used here, which is of order two.

\subsection{Propagation of periodic and regular nonlinear waves}\label{fenton}

In this test, we want to evaluate the dispersive properties of the model, along with the optimisation on the semi-discrete dispersion relation proposed in \S \ref{reldispdt}.

\medbreak

\noindent We consider the propagation of two-dimensional periodic and regular nonlinear waves, without change of form, over a flat bottom. For this situation, numerical reference solutions can be obtained by the so-called stream function method (see \cite{Fenton}). Unlike analytical wave theories (such as Stokes or cnoidal wave theories), this numerical approach is applicable whatever the shallowness and nonlinearity parameters $\mu$ and $\varepsilon$ are, and very accurate solutions can be obtained with a high number of terms in the Fourier series. These reference solutions are here obtained with the software Stream\_HT, implemented by Benoit \textit{et al.} \cite{MB}.

\medbreak

\noindent We consider a domain which covers one wave-length ($L = \lambda = 2\,m$), and a still water depth $h_0 = 1\,m$, so that the relative water depth is $kh_0 = \pi \approx 3.14$. The wave amplitude is $a = 0.01\,m$, so that the nonlinearity parameter is $a/h_0 = 0.01$. These conditions correspond to very dispersive and weakly nonlinear waves.

\medbreak

\noindent The domain is discretized with 50 cells ($\delta_x = 0.04\,$m), and a time-step $\delta_t = 0.03\,$s (corresponding to a Courant number $C_r=2.3$) is used during the simulations. Periodic conditions are imposed at the two lateral boundaries. The period computed with the stream function approach (at order 20) is $T = 1.133\,$s and the solutions obtained for the the water height and the averaged velocity are imposed as initial conditions. Numerical integration is performed over a duration of $25\,T$.

\medbreak

\noindent Two different values of $\alpha$ are considered: $\alpha=1.16$, corresponding to the local optimal value computed for $kh_0 = \pi$ as in \S \ref{sectdim}, and $\alpha=1.153$, corresponding to the local optimal value computed for $kh_0 = \pi$ and $\widetilde{\delta_t} = 0.094s$, as in \S \ref{reldispdt}. The local minimization approach is prefered since the considered wave is monochromatic.

\medbreak

\begin{figure}[!t]
\centering
\psfrag{x}[][][0.85]{}
\psfrag{h}[][][0.85]{$h$}
\psfrag{u}[][][0.85]{$u$}
\includegraphics[width=0.485\textwidth]{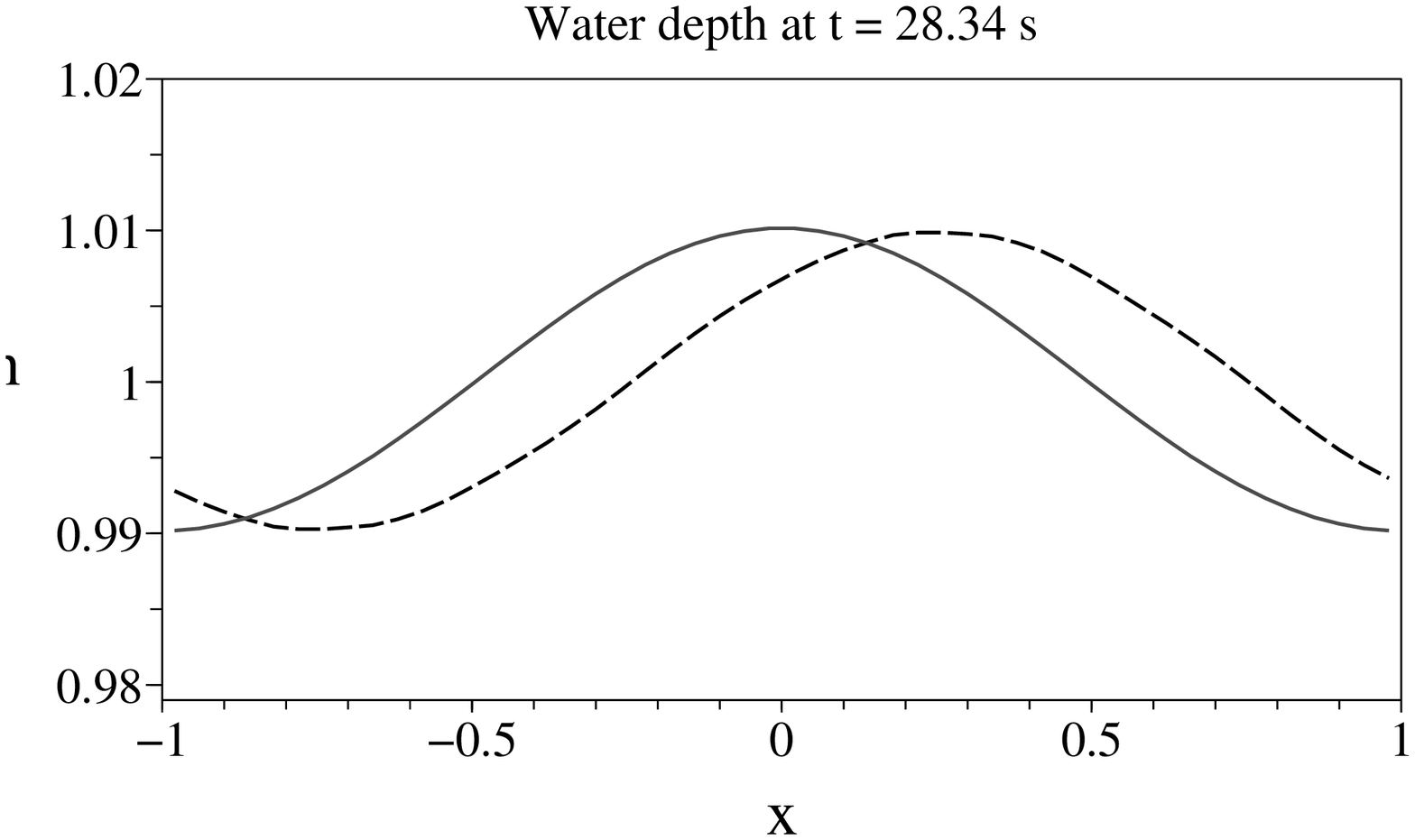}\;\;\includegraphics[width=0.485\textwidth]{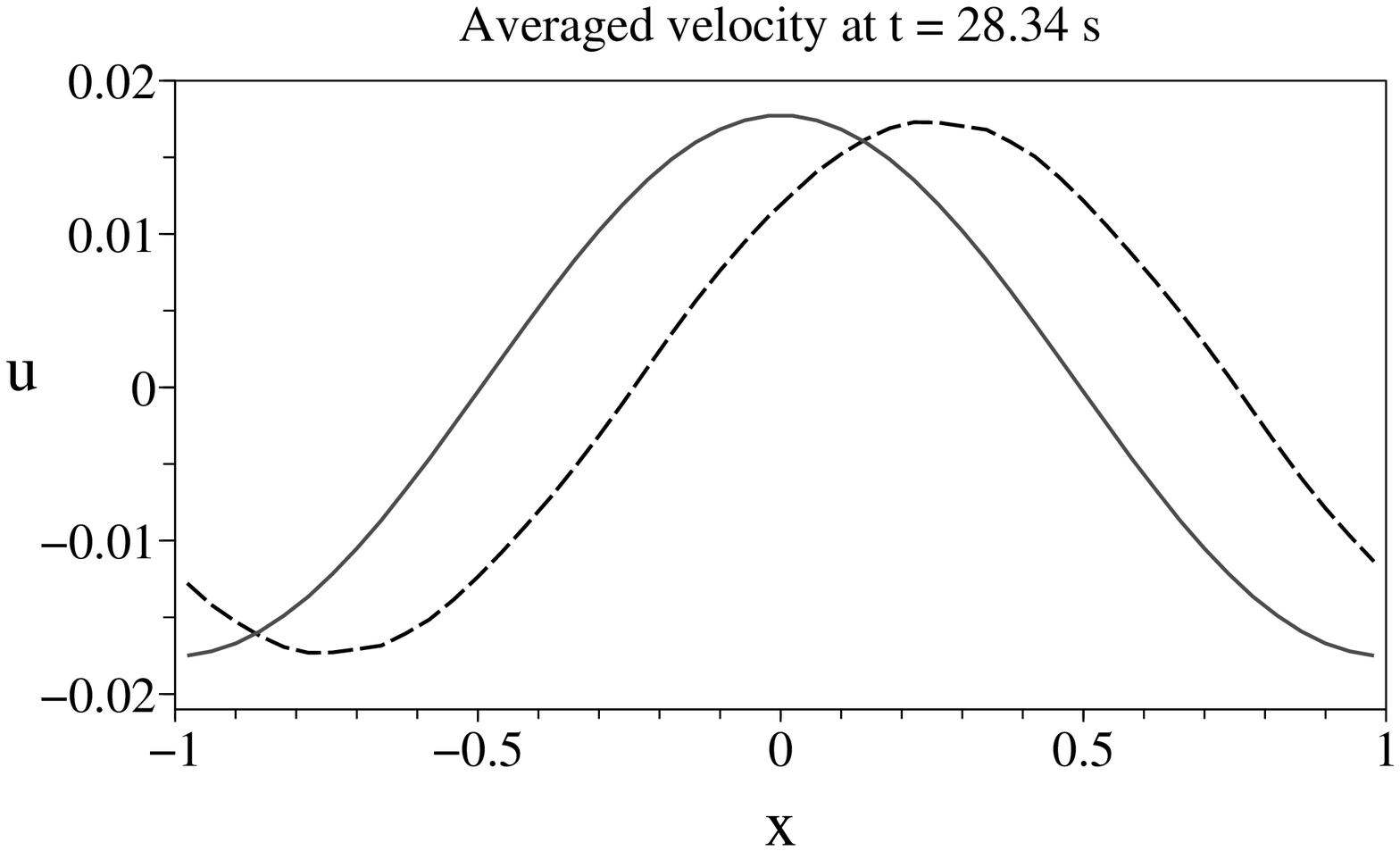}
\psfrag{x}[][][0.85]{$x$}
\psfrag{h}[][][0.85]{$h$}
\psfrag{u}[][][0.85]{$u$}
\includegraphics[width=0.485\textwidth]{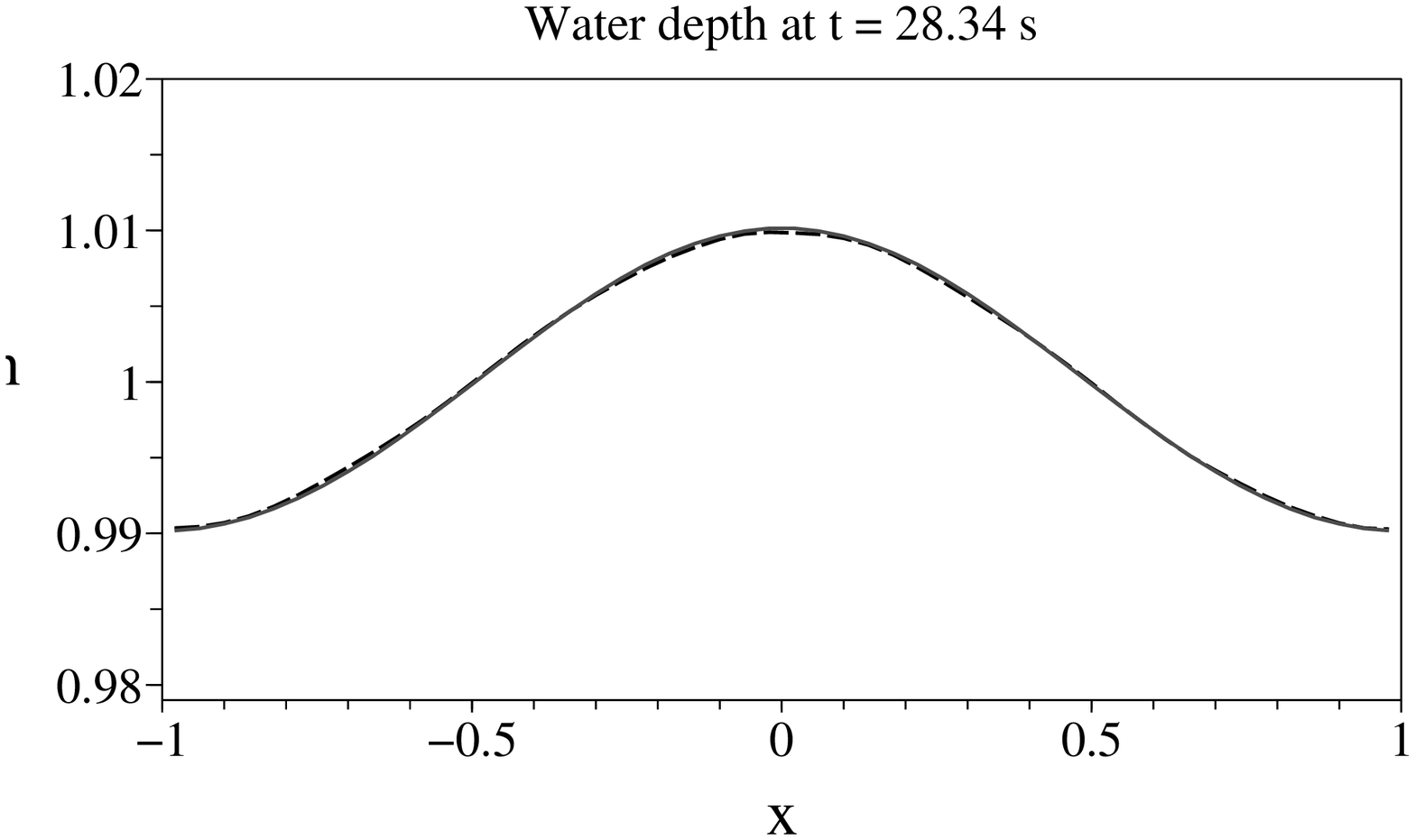}\;\;\includegraphics[width=0.485\textwidth]{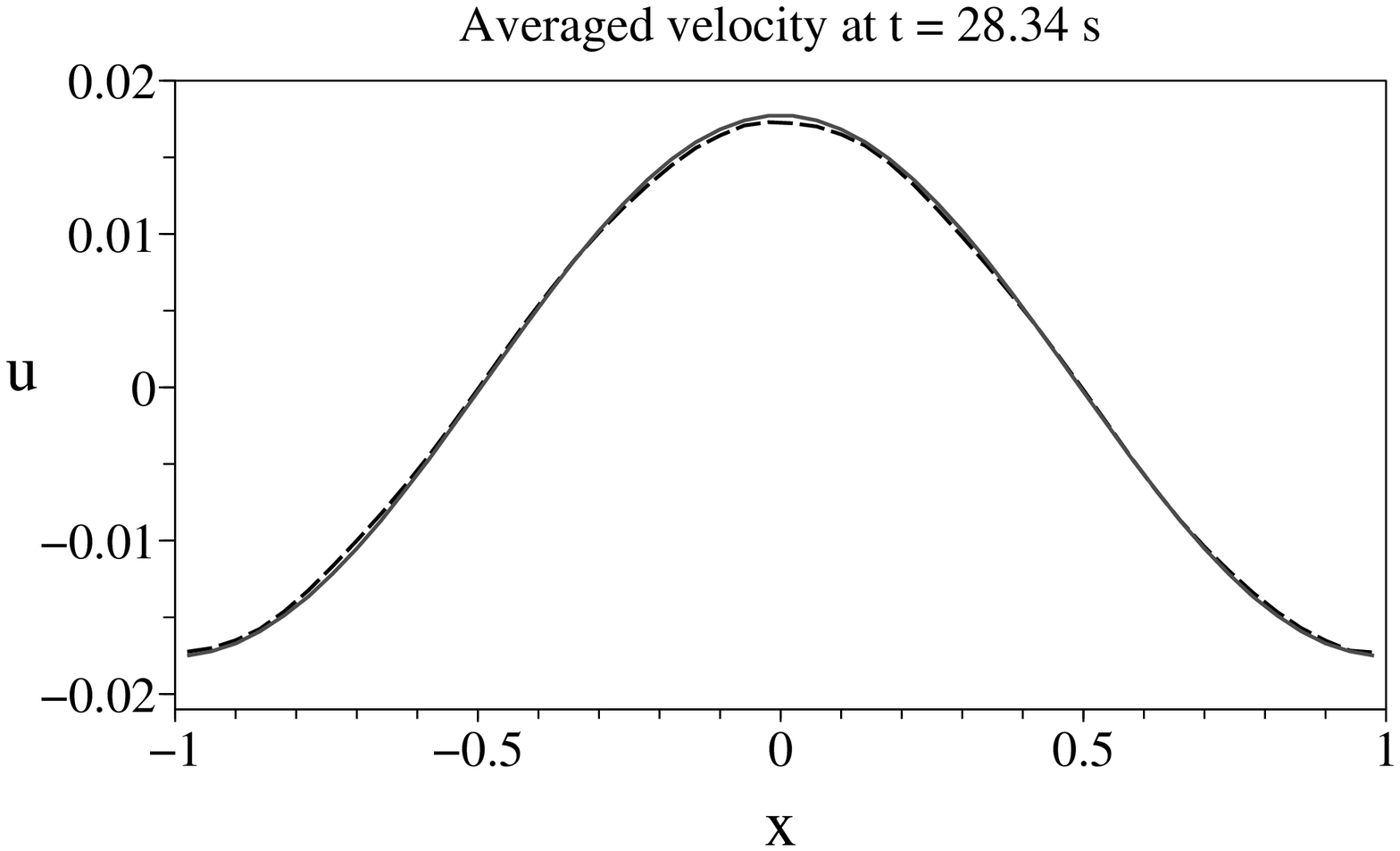}
\caption{Numerical results at $t=25\,T$, model with $\alpha=1.16$ (top) and $\alpha=1.153$ (bottom) in dashed line, reference solution in full line.}
\label{stream}
\end{figure}

Results after $25$ periods are gathered in Figure \ref{stream}, where the water height (left) and the averaged velocity (right) are plotted and compared to the reference solutions (propagating at constant speed and without change of form). The results obtained with $\alpha=1.153$ are seen to be in excellent agreement with the reference solution, whereas the ones obtained with $\alpha=1.16$ are less satisfying. However, we point out that these latter provide an overall good agreement with the reference solution, as shown in Table 1 where the relative error on the wave amplitude at $t=25\,T$ and the relative error on the wave celerity (using the phase shift between the solutions at $t=25\,T$) have been computed. To sum up, this test clearly demonstrates 1) the ability of our model to handle intermediate or even deep water waves (even for non-optimal values of $\alpha$), and 2) the interest of using the optimisation on the semi-discrete dispersion relation proposed in \S \ref{reldispdt} when dispersive waves are involved.

\begin{table}[!t]
\begin{center}
\begin{tabular}{|c|c|c|}
\hline
\vspace*{-0.9em} & & \\
 & Relative error on the & Relative error on the \\
 & wave amplitude & wave celerity \\
\hline
\vspace*{-0.9em} & & \\
$\alpha = 1.16$ & $1.8\;10^{-2}$ & $5.10^{-3}$ \\
\hline
\vspace*{-0.9em} & & \\
$\alpha = 1.153$ & $1.7.10^{-2}$ & $8.10^{-4}$ \\
\hline
\end{tabular}
\end{center}
\caption{Relative errors on the wave amplitude and the wave celerity between the model results and the reference solution (see Figure \ref{stream}) at $t = 25\,T$.}
\label{tab_erreur}
\end{table}

\subsection{Reflection of a solitary wave at a wall}\label{wallwb}

\begin{figure}[!t]
\centering
\includegraphics[width=\textwidth]{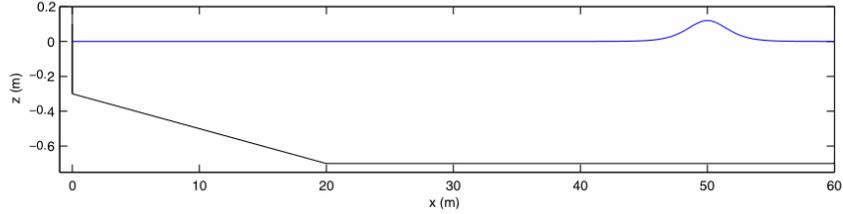}
\caption{Topography layout for the reflection of a solitary wave at a wall.}
\label{walltest}
\end{figure}

\begin{figure}[!t]
\includegraphics*[width=6cm, height=5cm, angle=0]{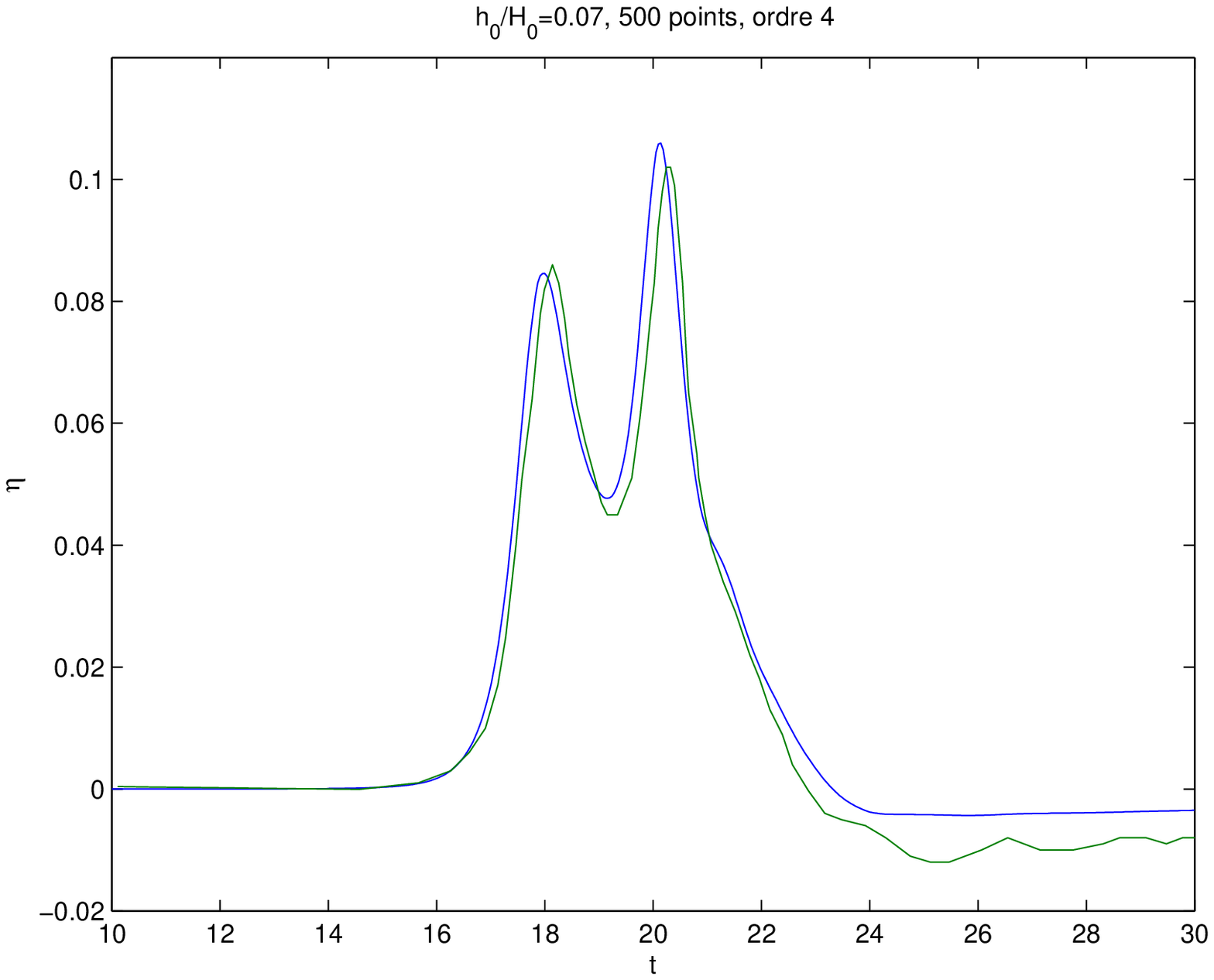}
\includegraphics*[width=6cm, height=5cm, angle=0]{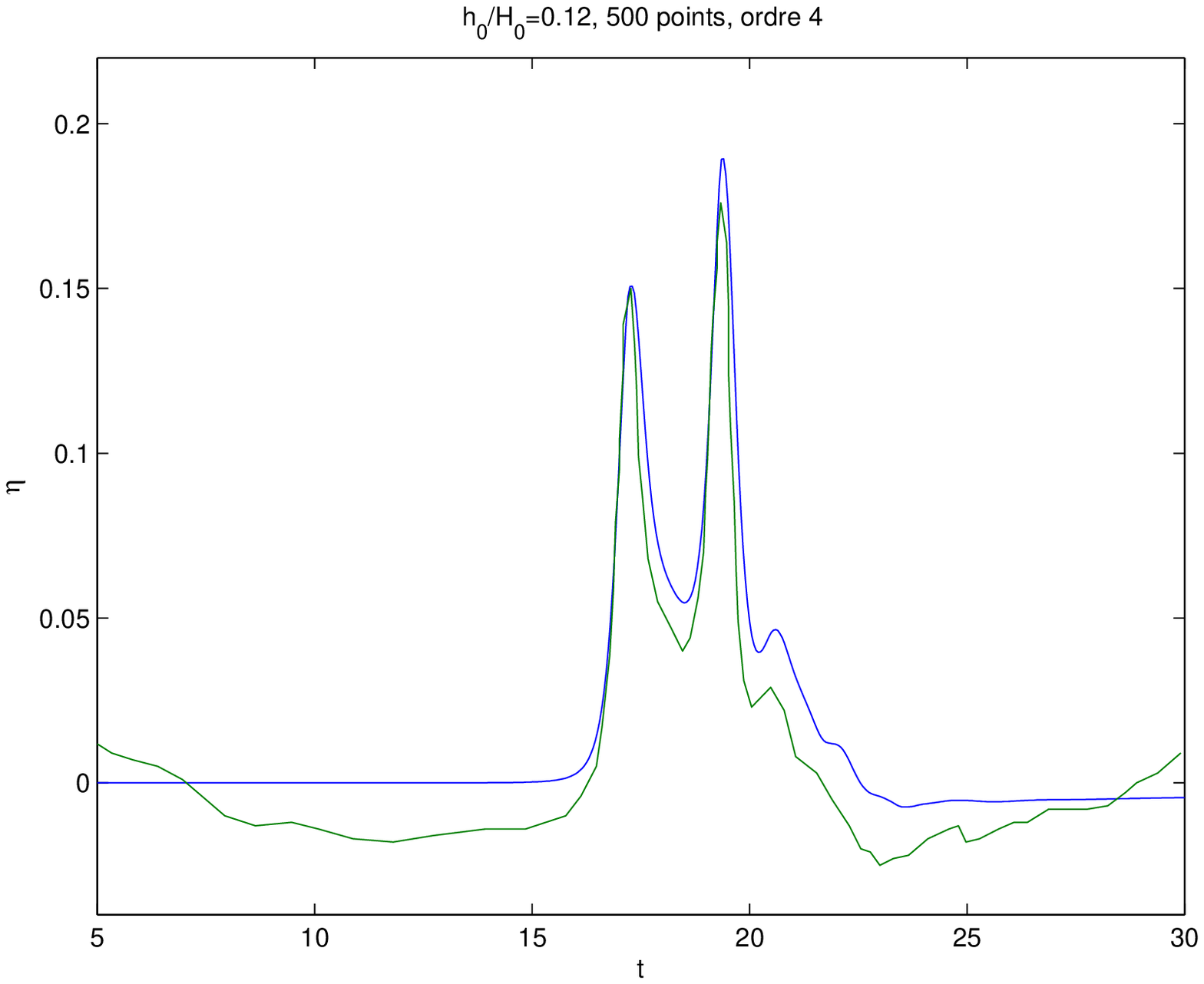}
\caption{Reflection of a solitary wave against a vertical wall. Time series of the free surface at $x=17.75\,m$ for $a=0.07$ (left) and $a=0.12$ (right): comparison between numerical results (blue line) and experimental data (green line).}
\label{wall}
\end{figure}

In this test, we compare numerical results with experimental data taken from \cite{wb}, for the propagation and reflexion of a solitary wave against a vertical wall. The depth profile, a sloping beach of 1:50  terminated by the wall located at $x=0\,m$, is depicted in Figure \ref{walltest}. The aim of this test is to study the full reflection of a non-breaking solitary wave propagating above a regular sloping beach, before reaching a vertical solid wall.

\medbreak

\noindent The spatial domain is $60\,m$ long, the initial solitary wave is centred at $x=50\,m$ and is propagating from right to left. The still water depth is $h_0 = 0.7 \,m$. The boundary condition at right is open, as there is no inflow. At the left boundary we use solid walls fully reflective conditions.

\medbreak

\noindent Two runs are performed with two different initial solitary wave amplitudes, given in terms of relative amplitude $a/h_0=0.1$ and $a/h_0=0.174$. The computational domain is discretized using 500 cells and a time step $\delta_t =0.05\,s$ is used.

\medbreak

Numerical results are shown as time series of the surface elevation, at a location near the solid wall ($x=17.75 \,m$). Experimental data are compared with numerical results on Figure \ref{wall}. We can observe the two expected peaks corresponding respectively to the incident and reflected waves. We can observe a very accurate matching between simulation and experimental data, even for the second simulation which involves a more complex non-linear propagation.

\subsection{Nonlinear shoaling of solitary waves propagating over a beach}\label{LEGI}

We investigate in this test the ability of the scheme to simulate the nonlinear shoaling of solitary waves over regular sloping beaches, which is a paramount in the study of nearshore propagating waves. This test is based on the experiments performed at the LEGI, in Grenoble (France) and reported in \cite{guibourg}. Solitary waves are generated in a $36\,m$ long wave-flume, following the procedure described in \cite{barthelemy_guizien}.

\medbreak

\noindent Free surface displacements were measured at various locations, using wave gau\-ges located just before breaking. Four solitary waves of different heights are generated (see Table \ref{tab_amplitude}), in order to account for various nonlinearity effects during propagation towards the shore.

\medbreak

\noindent All simulations are performed using $\delta_x = 0.025 \,m$ and  $\delta_t = 0.016 \,s$. The initial water depth is $h_0 = 0.25 \,m$ in the horizontal part of the channel.

\medbreak

\begin{table}[!t]
\begin{center}
\begin{tabular}{|m{5.2cm}|c|c|c|c|c|}
  \hline
  \multicolumn{6}{|c|}{Incident wave amplitude: $a_0/h_0 = 0.096$} \\
     \hline
     Gauge location (m) & 2.430 & 2.215 & 1.960 &1.740 & 1.502 \\
     (relative to the shoreline) & & & & & \\
     Relative amplitude error (\%) &  -1.6  & -2.5  & -5.5  & -7.1  &-10.9  \\
     \hline
   \multicolumn{6}{|c|}{ Incident wave amplitude: $a_0/h_0 = 0.298$ }\\
     \hline
      Gauge location (m) & 3.980 & 3.765 & 3.510 & 3.290 & 3.052\\
     (relative to the shoreline) & & & & & \\
      Relative amplitude error (\%) &   1.2  & -0.5  &  0.2 &  -0.2  &  0.04  \\
     \hline
   \multicolumn{6}{|c|}{ Incident wave amplitude: $a_0/h_0 = 0.456$ }\\
     \hline
      Gauge location (m) & 4.910 & 4.695 & 4.440 & 4.220 & 3.982\\
     (relative to the shoreline) & & & & & \\
      Relative amplitude error (\%)  &  3.6 & -0.3 & 1.1 & 0.5 & 2.2  \\
     \hline
   \multicolumn{6}{|c|}{ Incident wave amplitude: $a_0/h_0 = 0.534$ }\\ 
     \hline
     Gauge location (m) &  5.180 & 4.965 & 4.710 & 4.490 & 4.252\\
     (relative to the shoreline) & & & & & \\
     Relative amplitude error (\%) & 0.03 &   -0.1 &  -1.4 & -1.7 &  0.7 \\
     \hline
\end{tabular}
\end{center}
\caption{Location of wave gauges for solitary waves shoaling on a 1:30 sloped beach, and  relative error between the computed and measured wave  amplitudes at each gauge. \label{tab_amplitude}}
\end{table}
\begin{figure}[t]
\includegraphics*[width=\linewidth]{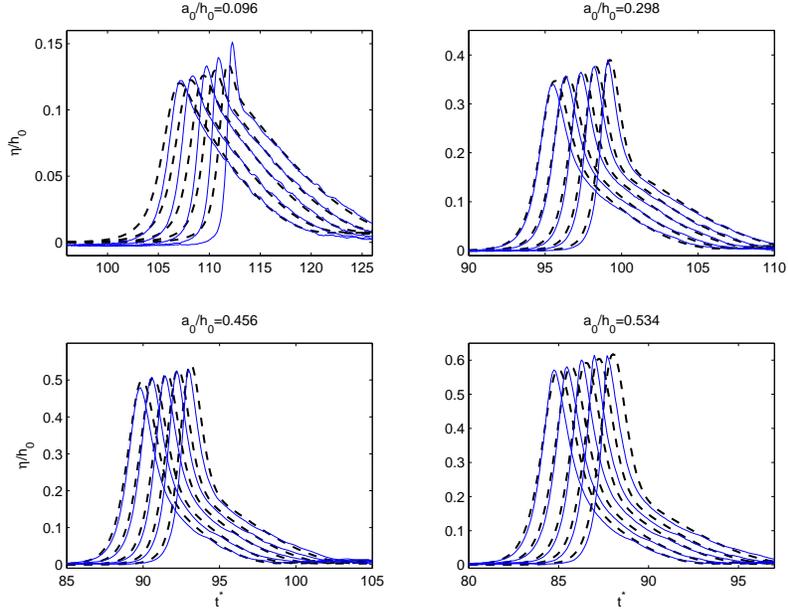}
\caption{Nonlinear shoaling of solitary waves propagating over a beach - Time series of the free surface elevation for the solitary wave propagating over the 1:30 sloping beach. (---) experimental data, (- - -) numerical results, with $t^*=t(g/h_0)^{1/2}$.}
\label{LEGI1}
\end{figure}

Results are shown  for each configuration in terms of time-series at the wave gauges locations in Figure \ref{LEGI1}. The relative error between computed and measured wave amplitudes is presented in Table \ref{tab_amplitude}.
The global agreement is good, both for the amplitude and shape of the solitary waves. Significant errors can be observed for the
less nonlinear case ($a_0/h_0 =  0.096$), but the discrepancies can  be partly explained by experimental problems, since it can
be observed that the water surface is not totally at rest before the propagation of the solitary wave .

\subsection{Run-up and run-down of a breaking solitary wave over a planar beach}\label{syno}

This test is based on experiments carried out by Synolakis \cite{ref_syno} for an incident solitary wave of relative amplitude $a_0/h_0 = 0.28$, which propagates and breaks over a planar beach with a slope of 1:19.85. Free surface elevations at different times are available thanks to video measurements.

\medbreak

\noindent The still water level in the horizontal part of the beach is $h_0=0.3 \,m$. The simulations are performed using the cell size $\delta_x = 0.08 \,m$ and  $\delta_t = 0.02 \,s$. Friction effects are important when the water becomes very shallow, as for instance in the run-up and run-down stage.  
To take into account this phenomenon, we introduce a quadratic friction term with the friction coefficient $f = 0.002$ (see Remark \ref{friction}).

\medbreak

\begin{figure}[!t]
\centering
\includegraphics[height=0.62\textheight]{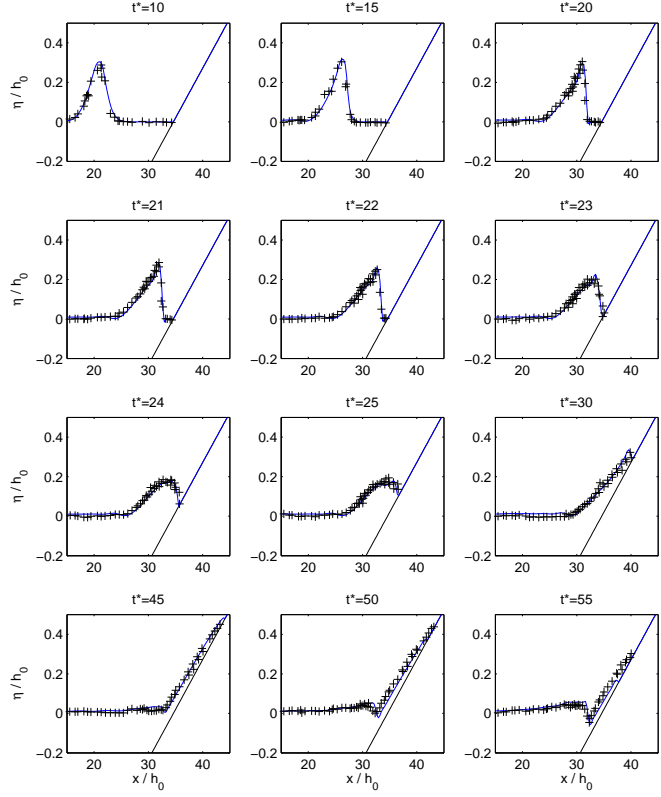}
\caption{Comparisons of model predictions (---) and experimental snapshots (+) for a breaking
solitary wave with non-dimensional initial incident amplitude $ a_0 /h_0 = 0.28$, on a $1:19.85$ constant slope
beach investigated by Synolakis (1987), with $t^*=t(g/h_0)^{1/2}$.}
\label{synolakis}
\end{figure}

The comparison between measured and computed waves is presented in Fi\-gure \ref{synolakis}. It shows a good agreement between model predictions and laboratory data and illustrates the ability of our model to reproduce shoaling, breaking, run-up and run-down, as well as the formation and breaking of the backwash bore, and this without any additional treatment.

\section{Conclusion}

In this work, the fully nonlinear and weakly dispersive Green-Naghdi model for shallow water waves is considered. The original model  is reformulated under a more convenient form, 
and variants with improved frequency dispersion depending on a parameter $\alpha$ are
derived.

A hybrid finite-volume and finite-difference method is then implemented, embedded in a splitting approach especially designed to describe the wide range of phenomena encountered in coastal oceanography. The first component of the model, regarded as a set of {\it hyperbolic} conservation laws, is discretized using an efficient and robust Godunov-like high-order accuracy finite-volume scheme, while high-order finite-differences are used for the {\it dispersive part} of the model.

The dispersive properties of the splitted semi-discrete scheme are carefully studied and are shown to approach the dispersion relation of the Green-Naghdi model at order two. Moreover, we use the explicit formula for the semi-discre\-tized dispersion relation to choose the best coefficient $\alpha$ for the frequency-improved GN model. This optimal value depends on the time step $\delta_t$, and is shown to provide better results than the standard choice based on the dispersion relation of the mathematical model. This is because this choice takes into account the dispersive effects due
to the time discretization.

In a last part, this new scheme is widely validated. Analytical solutions for propagating solitary waves are first considered and allow us to study the accuracy and convergence properties of this approach. In the following cases, numerical results are compared in an extensive way with both experimental data and reference solutions. In particular, the propagation and shoaling of highly nonlinear waves are successfully described, together with wave breaking and subsequent run-up and back-wash over a slopping beach. This clearly demonstrates the validity of this shock-capturing finite-volume based approach for dispersive waves, which appears therefore as a promising tool for the study of shallow water waves in coastal areas.

Following the steps of this study, next steps may concerns the derivation of dispersion optimized models \cite{proc}, the design of numerical sensor for the accurate detection of breaking waves and of course two-dimensional extension of the numerical scheme to study more realistic cases.

\section*{Acknowledgments}
This work has been supported by the ANR MathOcean, the project ECOS-CONYCIT action C07U01, and was also performed within the framework of the LEFE-IDAO program (Interactions et Dynamique de l'Atmosph\`ere et de l'Oc\'ean) sponsored by the CNRS/INSU. PhD thesis of M. Tissier is funded by the ANR MISEEVA.

\bigbreak


\begin{thebibliography}{00}

\bibitem{AL} B. Alvarez-Samaniego, D. Lannes, \textit{Large time existence for 3D water waves and asymptotics}, Invent. Math. 171(3) (2008) 485--541.
\bibitem{audusse} E. Audusse, F. Bouchut, M.-O. Bristeau, R. Klein, B. Perthame, \textit{A fast and stable well-balanced scheme with hydrostatic reconstruction for shallow water flows}, J. Comp. Phys. 25(6) (2004) 2050--2065.
\bibitem{berthon_marche} C. Berthon, F. Marche, \textit{A Positive Preserving High Order VFRoe Scheme for Shallow Water Equations: A Class of Relaxation Schemes}, SIAM J. Sci. Comput. 30(5) (2008) 2587--2612.
\bibitem{bar2004} E. Barthelemy, \textit{Nonlinear shallow water theories for coastal waves}, Surveys in Geophysics 25 (2004) 315--337.
\bibitem{MB} M. Benoit, M. Luck, C. Chevalier, M. B\'elorgey, \textit{Near-bottom kinematics of shoaling and breaking waves: experimental investigation and numerical prediction}, Proc. 28th Int. Conf. on Coastal Eng., Cardiff, UK (2002) 306--318.
\bibitem{bonn2007} P. Bonneton, \textit{Modelling of periodic wave transformation in the inner surf zone}, Ocean Engineering 34 (2007) 1459--1471.
\bibitem{bouchut} F. Bouchut, \textit{Nonlinear stability of Finite Volume Methods for Hyperbolic Conservation Laws
and Well-Balanced Schemes for Sources}, Frintiers in Mathematics, Birkhauser, 2004.
\bibitem{brocchini} M. Brocchini, I. Svendsen, R. Prasad, G. Bellotti, \textit{A comparison of two different types of shoreline boundary conditions}, Computer Methods in Applied Mechanics and Engineering 191(39-40) (2008) 4475--4496.
\bibitem{brocchini2} M. Brocchini, N. Dodd, \textit{Nonlinear shallow water equation modeling for coastal engineering}, J. Wtrwy. Port, Coast. and Oc. Engrg 134(2) (2008) 104--120.
\bibitem{camarri} S. Camarri, M.-V. Salvetti, B. Koobus, A. Dervieux, \textit{A low-diffusion MUSCL scheme for LES on unstructured grids}, Computer and Fluids 33(9) (2004) 1101--1129.
\bibitem{DL} F. Chazel, M. Benoit, A. Ern, S. Piperno, \textit{A double-layer Boussinesq-type model for highly nonlinear and dispersive waves}, Proc. R. Soc. Lond. A 465 (2009) 2319--2346.
\bibitem{proc} F. Chazel, D. Lannes, F. Marche, \textit{Numerical simulation of strongly nonlinear and dispersive waves using a Green-Naghdi model}, submitted to J. Sci. Comput. (2010).
\bibitem{Chen} Q. Chen, J. T. Kirby, R. A. Dalrymple, A. B. Kennedy, A. Chawla, \textit{Boussinesq modeling of wave 
transformation, breaking, and runup. II: 2d}, J. Wtrwy., Port, Coast., and Oc. Engrg. 126 (2000) 48--56.
\bibitem{cie2006} R. Cienfuegos, E. Barthelemy, P. Bonneton, \textit{A fourth-order compact finite volume scheme for fully nonlinear and weakly dispersive Boussinesq-type equations. Part I: Model development and analysis}, Int. J. Numer. Meth. Fluids 56 (2006) 1217--1253.
\bibitem{cie2007} R. Cienfuegos, E. Barthelemy, P. Bonneton, \textit{A fourth-order compact finite volume scheme for fully nonlinear and weakly dispersive Boussinesq-type equations. Part II: Boundary conditions and validations}, Int. J. Numer. Meth. Fluids 53 (2007) 1423--1455.
\bibitem{cie2009} R. Cienfuegos, E. Barthelemy, P. Bonneton, \textit{A wave-breaking model for Boussinesq-type equations including mass-induced effects}, J. Wtrwy., Port, Coast., and Oc. Engrg. 136 (2010) 10--26.
\bibitem{erduran} K. S. Erduran, S. Ilic, V. Kutija, \textit{Hybrid finite-volume finite-difference scheme for the solution of Boussinesq equations}, Int. J. Numer. Meth. Fluids. 49 (2005) 1213--1232.
\bibitem{godlewski_raviart} E. Godlewski, P.-A. Raviart, \textit{Numerical approximation of hyperbolic systems of conservation laws}, Applied Mathematical Sciences 118, Springer, 1996.
\bibitem{gre1976} A. E. Green, P. M. Naghdi, \textit{A derivation of equations for wave propagation in water of variable depth}, Journal of Fluid Mechanics 78(2) (1976) 237--246.
\bibitem{guibourg} S. Guibourg, \textit{Mod\'{e}lisation num\'{e}rique et exp\'{e}rimentale des houles bidimensionnelles en zone coti\`{e}re}, PhD Thesis, Universit\'{e} Joseph Fourier - Grenoble I, France, 1994.
\bibitem{barthelemy_guizien} K. Guizien, E. Barthelemy, \textit{Accuracy of solitary wave generation by a piston wave maker}, Journal of Hydraulic Research 40(3) (2002) 321--331.
\bibitem{Kennedy} A. B. Kennedy, Q. Chen, J. T. Kirby, R. A. Dalrymple, \textit{Boussinesq modeling of wave transformation, breaking, and runup. I: 1D}, J. Wtrwy., Port, Coast., and Oc. Engrg. 126 (1999) 39--47.
\bibitem{koba1989} N. Kobayashi, G. De Silva, K. Watson, \textit{Wave transformation and swash oscillation on gentle and steep slopes}, J. Geophys. Res. 94 (1989) 951--966.
\bibitem{Lannes-Bonneton} D. Lannes, P. Bonneton, \textit{Derivation of asymptotic two-dimensional time-dependent equations for surface water wave propagation}, Physics of Fluids 21(1) (2009) 016601.
\bibitem{lynett} P. Lynett, T. Wu, P. Liu, \textit{Modeling wave runup with depth-integrated equations}, Coastal Engineering 46 (2002) 89--107.
\bibitem{LGH} O. Le M\'etayer, S. Gavrilyuk, S. Hank, \textit{A numerical scheme for the Green-Naghdi model}, J. Comp. Phys. 229 (2010) 2034--2045. 
\bibitem{marche_bonneton} F. Marche, P. Bonneton, P. Fabrie, N. Seguin, \textit{Evaluation of well-balanced bore-capturing schemes for 2D wetting and drying processes}, Internat. J. Numer. Methods Fluids 53(5) (2007) 867--894.
\bibitem{mad1991} P. A. Madsen, R. Murray, O. R. Sorensen, \textit{A new form of the Boussinesq equations with improved linear dispersion characteristics}, Coastal Eng. 15 (1991) 371--388.
\bibitem{Madsen} P. A. Madsen, H. B. Bingham, H. Liu, \textit{A new Boussinesq method for fully nonlinear waves from shallow to deep water}, J. Fluid Mech. 462 (2002) 1--30.
\bibitem{Madsen2} P. A. Madsen, H. B. Bingham, H. A. Sch‰ffer, \textit{Boussinesq-type formulations for fully nonlinear and extremely dispersive water waves: derivation and analysis}, Proc. R. Soc. Lond. A 459 (2003) 1075--1104.
\bibitem{MS} J. Miles, R. Salmon, \textit{Weakly dispersive nonlinear gravity waves}, Journal of Fluid Mechanics 157 (1985) 519--531.
\bibitem{per1967} D. H. Peregrine, \textit{Long waves on a beach}, Journal of Fluid Mechanics 27 (1967) 815--827.
\bibitem{Fenton} M. M. Rienecker, J. D. Fenton, \textit{A Fourier approximation for steady water waves}, J. Fluid Mech. 104 (1981) 119--137.
\bibitem{sea1987} F. J. Seabra-Santos, D. P. Renouard, A. M. Temperville, \textit{Numerical and experimental study of the transformation of
a solitary wave over a shelf or isolated obstacle}, Journal of Fluid Mechanics 176 (1987) 117--134.
\bibitem{sug1969} C. H. Su, C. S. Gardner, \textit{Korteweg-de Vries equation and generalizations. III. Derivation
of the Korteweg-de Vries equation and Burgers equation}, J. Math. Phys. 10(3) (1969) 536--539.
\bibitem{ref_syno} C. E. Synolakis, \textit{The runup of solitary waves}, Journal of Fluid Mechanics 185 (1987) 523--545.
\bibitem{vanleer} B. van Leer, \textit{Towards the ultimate conservative difference scheme. V. A second-order sequel to Godunov's method}, J. Comput. Phys. 32 (1979) 101--136.
\bibitem{vanhof} B. Van't Hof, E. A. H Vollebregt, \textit{Modelling of wetting and drying of shallow water using artificial porosity}, Internat. J. Numer. Methods Fluids 48(11) (2005) 1199--1217.
\bibitem{wb} M. Walkley, M. Berzins, \textit{A finite element model for the two-dimensional extended Boussinesq equations}, Internat. J. Numer. Methods Fluids 39(10) (2002) 865--885.
\bibitem{wei1995} G. Wei, J. T. Kirby, S. T. Grilli, R. Subramanya, \textit{A fully nonlinear Boussinesq model for surface waves. Part 1. Highly nonlinear unsteady waves}, Journal of Fluid Mechanics 294 (1995) 71--92.
\bibitem{funwave} G. Wei, J. T. Kirby, \textit{A time-dependent numerical code for extended Boussinesq equations}, J. Wtrwy., Port, Coast., and Oc. Engrg. 120 (1995) 251--261.
\bibitem{wit1984} J. M. Witting, \textit{A unified model for the evolution of nonlinear water waves}, J. of Comput. Phys. 56(2) (1984) 203--236.
\bibitem{Zelt} J. A. Zelt, \textit{The run-up of nonbreaking and breaking solitary waves}, Coastal Engineering 15 (1991) 205--246.

\end{thebibliography}
\end{document}